\documentclass[11pt]{imsart}
\makeindex
\RequirePackage[numbers]{natbib}
\startlocaldefs
\usepackage{latexsym,amsmath,amsfonts}
\usepackage{mathrsfs}
\newcommand{\nref}[1]{(\ref{#1})}
\renewcommand{\textbf}[1]{\begingroup\bfseries\mathversion{bold}#1\endgroup}

\def\text#1{\hbox{#1}}

\def\endproof{\mbox{\ $\Box$}}

\newcommand{\E}{\mathbb{E}}
\def\1{\mbox{1\hspace{-.20em}I}}
\newtheorem{theorem}{Theorem}[section]
\newtheorem{proposition}{Proposition}[section]
\newtheorem{lemma}{Lemma}[section]

\newtheorem{remark}{Remark}[section]

\renewcommand{\textbf}[1]{\begingroup\bfseries\mathversion{bold}#1\endgroup}
\def\text#1{\hbox{#1}}

\def\1{\mbox{1\hspace{-.20em}I}}

\def\e{\varepsilon}
\def\r_\e{r_\epsilon}

\def\r{\right}

\def\bt{{\bf\theta}}

\def\bt{{\boldsymbol{\theta}}}

\def\Var{{\rm Var}}

\numberwithin{equation}{section}
\def\AArm{\fam0 \mathrm}%
\def\AAk#1#2{\setbox\AAbo=\hbox{#2}\AAdi=\wd\AAbo\kern#1\AAdi{}}%
\def\AAr#1#2#3{\setbox\AAbo=\hbox{#2}\AAdi=\ht\AAbo\raise#1\AAdi\hbox{#3}}%
\def\BBz{{\AArm Z\!\!\!Z}}%
\def \E{\hbox{\it I\hskip -2pt E}}
\def \P{\hbox{\it I\hskip -2pt P}}
\def \V{\mbox{\rm Var}}
\def \I{\hbox{\rm 1\hskip -3pt I}}

\def \ep{\epsilon}
\def \V{\mbox{\rm Var}}
\def \I{\hbox{\rm 1\hskip -3pt I}}
\endlocaldefs
\begin{document}
\begin{frontmatter}
% "Title of the Paper"
\title{Sharp detection of smooth signals in a high-dimensional sparse matrix with indirect observations}% \thanksref{t1}}
%\thankstext{t1}{This is an original survey paper}
\runtitle{Sharp detection of smooth signals in a sparse matrix }

\begin{aug}
\author{\fnms{Cristina} \snm{Butucea} %\thanksref{t1}
\ead[label=e1]{cristina.butucea@univ-mlv.fr}}
\address{Universit\'e Paris-Est Marne-la-Vall\'ee, \\
LAMA(UMR 8050), UPEMLV
F-77454, Marne-la-Vall\'ee, France\\ email:
cristina.butucea@univ-mlv.fr }

\and
\author{\fnms{Ghislaine} \snm{Gayraud}  %\thanksref{t2} \corref{Gayraud}
\ead[label=e2]{ghislaine.gayraud@utc.fr}}

\address{Universit\'e de Technologie de Compi\`egne $\&$ CREST, \\  BP 20529, 60205 Compi\`egne, France\\ email: ghislaine.gayraud@utc.fr}

%\thankstext{t1}{Cristina Butucea's research}

%\thankstext{t2}{Ghislaine Gayraud's research}

\runauthor{C. Butucea and G. Gayraud}

\affiliation{Universit\'e de Technologie de Compi\`egne $\&$ CREST \thanksref{t2}\\ Universit\'e de Marne la Vall\'ee \thanksref{t1}}

\end{aug}

% indicate corresponding author with \corref{}
% \author{\fnms{John} \snm{Smith}\thanksref{t1}\corref{Gayraud}\ead[label=e1]{smith@foo.com}\ead[label=e2,url]{www.foo.com}}
 %\thankstext{t1}{Thanks to somebody}
 %\address{line 1\\ line 2\\ \printead{e1}\\ \printead{e2}}

%\author{\fnms{Ghislaine} \snm{Gayraud} \corref{}\ead[label=e1]{ghislaine.gayraud@utc.fr}}
 % \thankstext{t2}{Ghislaine Gayraud's research was partially supported by the
%ANR-blanc "SP Bayes". \printead{e1}}
%\address{Universit\'e de Technologie de Compi\`egne and CREST, \\  BP 20529, 60205 Compi\`egne, France.}

%\and
%\author{\fnms{Yuri} \snm{Ingster} \thanksref{t2} \corref{}\ead[label=e2] {yurii_ingster@mail.ru}}
%\thankstext{t2}{\printead{e2}}
%\address{St. Petersburg State Electrotechnical University, 5, \\
%Prof. Popov str., 197376 St.Petersburg, Russia. \\
%\printead{e2}}

%\runauthor{???}

\begin{abstract}
We consider a matrix-valued Gaussian sequence model, that is, we observe a sequence of high-dimensional $M \times N$ matrices of heterogeneous Gaussian random variables $x_{ij,k}$ for $i \in\{1,...,M\}$, $j \in \{1,...,N\}$ and $k \in \mathbb{Z}$. The standard deviation of our observations is $\ep k^s$ for some $\ep >0$ and $s \geq 0$. 
 
We give sharp rates for the detection of a sparse submatrix of size $m \times n$ with active components. A component $(i,j)$ is said active if the sequence $\{x_{ij,k}\}_k$  have mean $\{ \theta_{ij,k}\}_k$ within a Sobolev ellipsoid of smoothness $\tau >0$ and total energy $\sum_k \theta^2_{ij,k} $ larger than some $r^2_\ep$. Our rates involve relationships between $m,\, n, \, M$ and $N$ tending to infinity such that $m/M$, $n/N$ and $\ep$ tend to 0, such that a test procedure that we construct has asymptotic minimax risk tending to 0.

We prove corresponding lower bounds under additional assumptions on the relative size of the submatrix in the large matrix of observations. Except for these additional conditions our rates are asymptotically sharp. Lower bounds for hypothesis testing problems mean that  no test procedure can distinguish between the null hypothesis (no signal) and the alternative, i.e. the minimax risk for testing tends to 1.
\end{abstract}



% history:
% \received{\smonth{1} \syear{0000}}

%\tableofcontents

\end{frontmatter}

\newpage

%%%%%%%%%%%%%%%%%%%%%%%%%%%%%%%%%%%%%%%
\section{Introduction}\label{Intro}
%%%%%%%%%%%%%%%%%%%%%%%%%%%%%%

Large matrices are used to model more and more applied problems in different areas such as signal theory, genomics, medical statistics. In case we observe large matrices of data on some period of time, we propose a procedure to test  whether a smaller submatrix only contains active components, that is smooth signal with some given smoothness and significant energy (measured by its $\mathbb{L}_2$-norm). This step should be taken as a preliminary step for dimension reduction.

This problem can be stated equivalently in the Gaussian sequence model of coefficients (say Fourier coefficients) of the signals. We propose to deal with the Gaussian sequence model, as it is easier for our computations and discuss later on the alternative interpretation as signal detection.
We include heterogeneous Gaussian observations in order to include the setup of indirect observations.

\bigskip
More precisely, we consider the following Gaussian sequence model 
\begin{eqnarray}\label{M1}
x_{ij,k}&=& \xi_{ij} \; \theta_{ij,k} + \epsilon \,\sigma_{ij,k} \, \eta_{ij,k},\; i \in I = \{1,...,M\}, \; j \in J = \{1,...,N\}, \; k \in \mathbb{Z},
\end{eqnarray} 
where $\{\eta_{ij,k}\}_{i \in I, j \in J, k \in \BBz}$ is a sequence of independent standard Gaussian random variables, $\sigma_{ij,k} >0$ are known and $\epsilon>0$ is the noise level. 
The $ M \times N $-matrix  $\xi = [\xi_{ij}]_{(i,j) \in I \times J},$ is deterministic  
(unknown) and has elements in $\{0,1\}$.  

In what follows,  the standard deviations $\sigma_{ij,k}$ 
are supposed to be the same for all components of the matrix, that is
$\sigma_{ij,k} = \sigma_k$ for all $k$ do not depend on $(i,j)$ in $I\times J$. 
We assume throughout the paper that, for some fixed given $s\geq 0$,
\begin{equation*}\label{noise}
\sigma_k \sim \vert k \vert^{s}, \mbox{ for large enough integer values of } \vert k \vert.
\end{equation*}
On the one hand, the case $s=0$ reduces to the case of direct observations of the signal. In that case, we could generalize our results to unknown (but constant) variance $\sigma$.
On the other hand, the case $s>0$ corresponds to signals observed in inverse problems like convolution with some independent noise, tomography etc. 

The polynomial behaviour of $\sigma_k$ as $k$ grows to infinity corresponds to mildly ill-posed inverse problems. We refer to \cite{CGPT} for more discussion on the relation between the sequence model with increasing variance and inverse problems in the Gaussian white noise model.

\bigskip

The matrix-valued sequence 
$\overline{\boldsymbol{\theta}} = [\xi_{ij}\{\theta_{ij,k} \}_{k \in \BBz} ]_{(i,j) \in I \times J}$
is the quantity of interest. We want to detect from observations in
the model (\ref{M1}) whether there is only noise or whether there are 'active components' in $\overline{\boldsymbol{\theta}}$, corresponding to $(i,j)$ where  $\xi_{ij}=1$. When a component $(i,j)$ is active, we assume that the corresponding sequence $\{\theta_{ij,k}\}_k$ belongs to a Sobolev ellipsoid and has significant total energy, i.e., $\{\theta_{ij,k}\}_k \in \Sigma (\tau, r_\epsilon)$, $\tau >0,\, r_\epsilon >0$,
where 
\begin{equation} \label{Sobolev.class}
\Sigma (\tau, r_\epsilon) =\{ \boldsymbol{\theta} \in l_2 (\BBz) : \; (2 \pi)^{2 \tau}  \sum_{k\in \BBz}| k|^{2\tau}  \;  \theta_{k}^2 \leq 1;   \sum_{k\in \BBz}
 \; \theta_{k}^2 \geq r_{\epsilon}^2
 \}
\end{equation}

In this paper, we assume that $\xi$ has a specific structure, i.e., it belongs to  
  \begin{eqnarray*}
  T_{M,N}(m,n) & = & \left\{ \xi \mbox{ matrix of size } M \times N : \exists \, A_\xi \subseteq I,\, \# A_\xi=m \mbox{ and } \exists \, B_\xi \subseteq J, \, \# B_\xi=n \right.\\
  %\left\{ \xi \in S_{M,N}(m,n): \exists \, A \subseteq I,\, \# A=m \mbox{ and } \exists \, B \subseteq J, \, \# B=n \right.\\
&& \left.\mbox{ such that } \xi_{ij} = \1((i,j) \in A_\xi\times B_\xi) \right\},
\end{eqnarray*}
where the non null elements form a submatrix with $m$ rows and $n$ columns. We shall always denote by $A_\xi$ and $B_\xi$ those rows and columns where the matrix $\xi \in T_{M,N}(m,n)$ has non null elements.

% {\bf Maybe} We may also consider matrices $\xi$  belonging to
% \begin{eqnarray*}
%   U_{M,N}(m,n) & = & \left\{ \xi \in S_{M,N}(m,n): \exists \, k \in \{1,...,M-m\}  \mbox{ and } \exists \, j \in \{1,...,N-n\}, \,  \right.\\
% && \left.\mbox{ such that } \xi_{ij} = I((i,j) \in \{k,...,k+m\}\times \{j,...,j+n\} \right\},
%   \end{eqnarray*}
% that is such that the non null elements form a block-submatrix with adjacent $m$ rows and $n$ columns.

The testing problem  of interest is
the following
\begin{eqnarray*}
H_0 &\; : \; & \overline{\boldsymbol{\theta}}  = 0 \;\;\;\;\;\;\\
H_1(\tau,r_\epsilon) &\; : \; & \overline{\boldsymbol{\theta}}  \in \Theta_{M,N}(\tau,r_\epsilon,m,n),
\end{eqnarray*}
where, for $\tau, \, r_\epsilon > 0$  and for $m,\,n,\, M$ and $N$ large, such that $m \leq M$ and $n \leq N$, we define
\begin{eqnarray*} 
\Theta_{M,N}(\tau,r_\epsilon,m,n)&=& \{ \overline{\boldsymbol{\theta}}=[\xi_{ij} \{\theta_{ij,k}\}_{k \in \BBz}]_{(i,j)\in I\times J} : \; \xi \in T_{M,N}(m,n), \\
&& \;
  \mbox{{\rm and for all} } (i,j) \in A_\xi \times B_\xi, \; 
\{\theta_{ij,k}\}_k \in \Sigma (\tau, r_\epsilon)
  \}.
\end{eqnarray*}
The alternative hypothesis consists of matrices of size $M \times N$ containing mainly noise, except for elements in some submatrix of size $m \times n$ containing sequences of Fourier coefficients of signals with Sobolev smoothness $\tau$ and energy ($\mathbb{L}_2$ norm) significantly large (larger than $r_\ep$).

\begin{remark} 
We may also assume that the matrix $\xi $ has entries either 0 or 1, such that $\sum_{(i,j) \in I\times J} \xi_{ij}=m\times n$. That means that we know the number of non null elements of 
the matrix $\xi$ but they can be found anywhere in the matrix. This case is exactly the vector case previously studied by \cite{GI} under the sparsity condition that  the number of active components $mn$ satisfies  $mn=(MN)^{1-b}$, where $b \in (0,1)$ corresponds to the sparsity index. 
%  {\bf give sparsity assumptions in GI and more precise results after the test is defined}
\end{remark}

\bigskip

Denote by $\P_0$ and    $\P_{\overline{\boldsymbol{\theta}}}$ the distributions  under  the null and  the alternative, respectively.
Denote also by $\E_0$, $\V_0$ and  $\E_{\overline{\boldsymbol{\theta}}}$,  $\V_{\overline{\boldsymbol{\theta}}}$ the expected values and variances  with respect to $\P_0$ and $\P_{\overline{\boldsymbol{\theta}}}$, respectively.
Set $\boldsymbol{\theta}_{ij} = \{\theta_{ij,k}\}_{k \in \BBz}$; 
%The notation  $\P_{{ \boldsymbol{\theta}_{ij}}}$, $\E_{{ \boldsymbol{\theta}_{ij}}}$ and $\V_{{ \boldsymbol{\theta}_{ij}}}$ will also be used:  they are related  to the joint distribution of the  observations $\{x_{ij,k}\}_{k \in \BBz}$ when $\xi_{ij}=1$.  
indices of probabilities, expectations or variances  which are expressed in terms of non-overlined subsequences of  ${\boldsymbol \theta}$ mean that they   correspond to active components.

For any test procedure $\psi$, that is, any measurable function with respect to the observations, taking values in $[0,1]$, set
$ \omega(\psi)=\E_0( \psi)$ its type I error probability and  $ \beta (\psi,\Theta_{M,N}(\tau,r_\epsilon,m,n))=\displaystyle{\sup_{{\overline{\boldsymbol{\theta}}} \in \Theta_{M,N}(\tau,r_\epsilon,m,n) }}
\E_{\overline{\boldsymbol{\theta}}} (1 - \psi)$ its maximal type II error probability over the set  $\Theta_{M,N}(\tau,r_\epsilon,m,n)$. 
Let us denote by 
$$ 
\gamma (\psi,  \Theta_{M,N}(\tau,r_\epsilon,m,n))=\omega(\psi) + \beta (\psi, \Theta_{M,N}(\tau,r_\epsilon,m,n))
$$
the total error probability of $\psi$  and denote by  $\gamma$
 the minimax total error probability over  $\Theta_{M,N}(\tau,r_\epsilon,m,n)$  which is defined by
\begin{eqnarray*} \gamma &:= & \gamma (\Theta_{M,N}(\tau,r_\epsilon,m,n))=
 \displaystyle{\inf_{\psi}} \gamma(\psi,\Theta_{M,N}(\tau,r_\epsilon,m,n)), \label{gamma} \end{eqnarray*}
where the infimum     is taken over all  test procedures.
We can not distinguish $H_0$ and $H_1(\tau,r_\epsilon) $ if
$ \gamma \rightarrow 1$ and
distinguishability occurs if there exists $\psi$ such that
 $\gamma (\psi,  \Theta_{M,N}(\tau,r_\epsilon,m,n)) \rightarrow 0$.

The aim of this paper is to derive distinguishability conditions and separation rates for alternatives $\Theta_{M,N}(\tau,r_\epsilon,m,n)$ and to determine statistical procedures $\psi$ (at least of asymptotic $\alpha$-level) which achieve these separation rates.  By separation rates, we mean
a family $\tilde{r_\ep}$ such that \\
$\left\{ \begin{array}{lll}
 \gamma \rightarrow 1  & \mbox{{\rm if }} & \displaystyle{\frac{r_\epsilon}{\tilde{r_\ep}}} \rightarrow 0, \\
 && \\

\gamma (\psi, \Theta_{M,N}(\tau,r_\epsilon,m,n)) \rightarrow 0 & \mbox{{\rm if }} & \displaystyle{\frac{r_\epsilon}{\tilde{r_\ep}}} \rightarrow +\infty. \end{array} \right.$\\

By sharp separation rates, we mean
a family $\tilde{r_\ep}$ such that \\
$\left\{ \begin{array}{lll}
  \gamma \rightarrow 1  & \mbox{{\rm if }} & \displaystyle{\limsup \frac{r_\epsilon}{\tilde{r_\ep}}}<1, \\
 && \\

\gamma (\psi, \Theta_{M,N}(\tau,r_\epsilon,m,n)) \rightarrow 0 & \mbox{{\rm if }} & \displaystyle{\liminf \frac{r_\epsilon}{\tilde{r_\ep}}} > 1. \end{array} \right.$\\

The  asymptotics for model (\ref{M1}) are  given by   $\epsilon \rightarrow 0$ and, as we are mainly  interested
in high-dimensional settings, by 
\begin{eqnarray}
m,\, n,\, M \mbox{ and } N \rightarrow +\infty , \quad
p=\frac{m}{M} \rightarrow 0, \quad q=\frac{n}{N} \rightarrow 0. \label{Sp-C}
\end{eqnarray}

 Here and later asymptotics and  symbols  $o$, $O$, $\sim$ and $\asymp$ are considered under $\epsilon \rightarrow 0$ and $m,\, n, \, M $ and $N$ such that $(\ref{Sp-C})$ holds.
 %, except when
%it is explicitly specified like, for e.g.,  $o_{N,M}$ which holds under $N \rightarrow +\infty $ and $M \rightarrow +\infty$ and %$o_\ep$ when $\ep \to 0$. \\

The plan of the paper is as follows. Section~\ref{sec:motivation} explains how this model is related to the multivariate Gaussian white noise model and how the inverse problem reduces to heterogenous observations in our Gaussian sequence model. In Section~\ref{sec:upper.bounds} we define the test procedure and give sufficient conditions such that the minimax risk for testing tends to 0. The construction of our test procedure involves solving an optimization problem. Section~\ref{sec:lower.bounds} presents the lower bounds for our problem and proofs are given in Section~\ref{sec:Proofs} and the Appendix.

%%%%%%%%%%%%%%%%%%%%%%%%%%%%%%%%%
\section{Sparse high-dimensional signal detection}\label{sec:motivation}
%%%%%%%%%%%%%%%%%%%%%%%%%%%%%%%%%

Let us see that the  previous problem arises  in some classical statistical models and hence,  it has a different interpretation. When dealing with high-dimensional data, we model functions of many variables with additive models. For many situations where additive models are employed see Stone~\cite{St.85} and references therein. Let us consider the multivariate Gaussian white noise model
\begin{equation}\label{gwn}
dX(t) = f(t) dt + \epsilon \cdot dW(t),\quad t \in [0,1]^d, \; d \in \mathbb{N},
\end{equation} 
$\e >0$ and $W(t)$ is the Wiener process. When estimating $f$ in a nonparametric model, the curse of dimensionality makes the rates exponentially slow for large dimension $d$. Additive models, where $f(t) = \sum_{j=1}^d f_j(t_j)$, $t_j \in [0,1]$ and $\int_0^1 f_{j} = 0$ for all $j$ from 1 to $d$, are estimated with much faster rates, but the global estimation risk still grows in a linear way with $d$. 
It is assumed in \cite{GI} that the univariate signal functions $f_j$ belong to a class $\mathcal{S}(\tau,r_\ep)$, i.e., it has Sobolev smoothness $\tau$ and total energy $\int_0^1 |f_j|^2$ larger than $r_\ep^2$. A function $f$ is Sobolev smooth if it belongs to $\mathbb{L}_2 ([0,1])$ such that $\int |\tilde f(u)|^2 (2\pi |u|)^{2\tau} du \leq 1$ (where $\tilde f$ is the caracteristic function of a function $f$) and $\tau$ is called its smoothness.

If we need to cope with very high dimension $d$, sparsity assumptions help reduce the dimension.
In Gayraud and Ingster \cite{GI}, it was assumed that only $d^{1-b}$ for some $0<b<1$  coordinates are significantly active, i.e., 
$f(t) = \sum_{j=1}^d \xi_j f_j(t_j)$, $\xi_j \in \{0,1\}$ for all $j$ from 1 to $d$ such that $\sum_j \xi_j = d^{1-b}$. They solved the following test problem: 
\begin{eqnarray*}
H_0 &\; :  & \mbox{all} \;\xi_{j} = 0 , \, \mbox{ (no signal is detected in data)}\;\;\;\;\;\;\\
H_1(\tau,r_\epsilon) &\; : & \mbox{there exists } d^{1-b} \mbox{ values of $j$ where } \xi_j=1 \mbox{ and } f_j \in \mathcal{S} (\tau,r_\ep).
\end{eqnarray*}
Different sharp detection rates were obtained along the values of $0<b<1$.

\bigskip

In our paper, we assume a sparse matrix structure for our additive model:
\begin{equation}\label{spadd}
f(t) = \sum_{i=1}^M \sum_{j=1}^N \xi_{ij} f_{ij} (t_{ij}),\quad t_{ij} \in [0,1]
\mbox{ and } \xi \in T_{M,N}(m,n),
\end{equation}
such that $\int_0^1 f_{ij} = 0$ for all $i$, $j$.
We call the component $(i,j)$ active if $\xi_{ij}=1$ and, in that case, we suppose that the signal in that coordinate belongs to the class $\mathcal{S}(\tau,r_\ep)$.

Let us reduce the sparse additive model \nref{gwn} such that \nref{spadd} holds to our initial model. Consider
$\{\varphi_k\}_{k\in\mathbb{Z}}$ an orthonormal basis of $\mathbb{L}_2[0,1]$ such that $\varphi_0 \equiv 1$ (e.g., the Fourier basis). Define the multivariate orthonormal family, for $t\in[0,1]^{M\times N}$,
$$
\Phi_{ij,k}(t) = \varphi_k(t_{ij}) \cdot \prod_{(l,h) \not = (i,j)} \varphi_0(t_{lh}) = \varphi_k(t_{ij}).
$$
Then, project the signal in \nref{gwn} on these functions:
\begin{eqnarray*}
x_{ij,k}&:= & \int_{[0,1]^{M\times N}} \Phi_{ij,k}(t) dX(t) \\
&=&\int_{[0,1]^{M\times N}} \Phi_{ij,k}(t) f(t)dt +\ep \cdot \int_{[0,1]^{M\times N}} \Phi_{ij,k}(t) dW(t)\\
&=& \xi_{ij} \int_0^1 \varphi_k(t_{ij}) f_{ij}(t_{ij}) dt_{ij} + \ep \cdot \eta_{ij,k},
\end{eqnarray*}
where $\{\eta_{ij,k}\}$ are i.i.d.  standard Gaussian random variables. We get our initial model for $\theta_{ij,k} = \int_0^1 \varphi_k f_{ij}$ and $\sigma_{k} \equiv 1$. 

Therefore, our test problem can be written:
\begin{eqnarray*}
H_0 &\; :  & \mbox{all} \;\xi_{ij} = 0 \mbox{ (no signal is detected in data)}\;\;\;\;\;\;\\
H_1(\tau,r_\epsilon) &\; : & \mbox{there exists } \xi \in T_{M,N}(m,n) \mbox{ and for } \xi_{ij}=1 \mbox{ it holds that } f_{ij} \in \mathcal{S}(\tau,r_\ep),
\end{eqnarray*}
i.e., there exists a matrix $\xi$ in $T_{M,N}(m,n)$ such that the signal in active coordinates $(i,j)$ has Sobolev smoothness $\tau$ and total energy larger than $r_\ep^2$.

\bigskip

The variance of our observations are allowed to increase $\sigma_k \sim \vert k \vert^s$, $s\geq 0$. Indeed, let us suppose that our additive model is observed as an inverse problem. That means that we observe
\begin{equation}\label{invgwn}
dX(t) =  Kf (t) dt +\ep \cdot dW(t), \quad t=[t_{ij}]_{i,j} \in [0,1]^{M\times N}
\end{equation}
for some linear operator $K$, with $f$ given as in \nref{spadd} and such that $\int_0^1 Kf_{ij}=0$. In the convolution model, for example, the signal is observed with an additive independent noise having density $g$, than $Kf(y) = \int f(y-u) g(u) du$.

We suppose that $K^* K$ is a compact operator having eigenvalues $\sigma_k^{-2}$ decreasing polynomially to 0 as $k$ tends to infinity. This corresponds to mildly ill-posed inverse problems. Whereas, in the case of well-posed inverse problems, $\sigma_k^2 \leq \sigma^2$ form a bounded sequence.

Then, we consider a singular value decomposition of $K$, that is families of orthonormal functions $\{\varphi_k \}_k$ and $\{\psi_k \}_k$ such that $K \varphi_k = \sigma_k^{-1} \psi_k$ and $K \psi_k = \sigma_k^{-1} \varphi_k$. Therefore, let $\psi_k\equiv 1$ and $\Psi_{ij,k}(t) = \psi_k(t_{ij})$,
and project \nref{invgwn} on this family:
\begin{eqnarray*}
y_{ij,k}&:=& \sum_{l=1}^M \sum_{h=1}^N \xi_{lh}\int_{[0,1]^{M\times N}} \Psi_{ij,k}(t)  Kf_{lh} (t_{lh}) dt_{lh} +\ep \cdot \int_{[0,1]^{M\times N}} \Psi_{ij,k}(t) dW(t)\\
&=& \xi_{ij} \int_0^1 \psi_{k}(u) Kf_{ij}(u) du + \ep \cdot \eta_{ij,k}.
\end{eqnarray*}
Note, moreover, that $\int_0^1 \psi_{k} \cdot Kf_{ij} = \int_0^1 K^*\psi_{k} \cdot f_{ij} =  \sigma_k^{-1} \int_0^1\varphi_{k} \cdot f_{ij} = \sigma_k^{-1} \theta_{ij,k}$. Then, let $x_{ij,k} = \sigma_k y_{ij,k}$ to get the model \nref{M1}.

Note that Butucea and Ingster \cite{BI} studied the particular case where $\theta_{ij,k} = a \1(k=0)$  and the variance of the noise is a given fixed $\sigma$. If we have in mind the Fourier basis, it comes down to studying periodic signals. The asymptotic rates for testing were given in terms of $n,\, m,\, N$ and $M$. Here, we replace  the periodic signal with arbitrary smooth signal. Moreover, we add here the case of heterogeneous variables which include mildly ill-posed inverse problems.

%%%%%%%%%%%%%%%%%%%%%%%%%%%%%%%%%%%%%%%%%%%%%%%%%%%%%%
\section{Testing procedures and their asymptotic behaviour}\label{sec:upper.bounds}
%%%%%%%%%%%%%%%%%%%%%%%%%%%%%%%%%%%%%%%%%%%%%%%%%%%%%%

Consider the following  family of weighted $\chi^2$-type statistics: for $(i,j)$ in $I \times J$
\begin{eqnarray*}
t_{ij,w} &=& \sum_{k\in \mathbb{Z}} w_k \left( (\frac{x_{ij,k}}{\epsilon \sigma_{k}})^2-1 \right), \label{tj}
%& =&  \sum_{k\in \BBz} ,
\end{eqnarray*}
where $(w_k)_k$ is a sequence of weights such that $w_k\geq 0$ forall $k \in \mathbb{Z}$ and $\sum_{k \in \mathbb{Z}} w_k^2 = 1/ 2$.   

In order to define the weights $\{w_k^\star\}_{k \in \mathbb{Z}}$ that will appear in the optimal test procedure, we have to solve the following extremal problem. Recall that  $\Sigma(\tau, r_\epsilon)$ denotes the Sobolev ellipsoid defined in  (\ref{Sobolev.class}), with $\tau >0$ and $r_\epsilon >0$, and $\{\sigma_k\}_{k \in \mathbb{Z}}$ is a sequence of positive real numbers. 
We define the sequences $\{w_k^\star\}_{k \in \mathbb{Z}}$ and $\{ \theta^\star_k\}_{k \in \mathbb{Z}}$ as solutions to the following optimization program:
%\mathop{\sum_{i \in I_k, \; j \in I_k}}_{ i
%\neq j }
\begin{eqnarray} \label{max-min_pb}
\sum_{k \in \mathbb{Z}} w_k^\star \left (\frac{\theta_k^\star}{\ep\sigma_k} \right)^2
 = \sup_{\left\{\begin{array}{ll} (w_k)_{k } \in l^2 (\mathbb{Z}): & w_k \geq 0; \\
 &  \sum_k w_k^2=\frac 1 2\end{array}\right\} } \quad \inf_{\{\theta_k\}_k \in \Sigma (\tau,r_\ep)}
\sum_{k \in \mathbb{Z}} w_k \left( \frac{\theta_k}{\ep \sigma_k} \right)^2 .
\end{eqnarray}
Let us denote by $V_\epsilon:=V_\epsilon(r_\epsilon) = \sum_{k \in \mathbb{Z}} w_k^\star \left ( \theta_k^\star / \sigma_k \right)^2$, so that $a(r_\epsilon): = V_\epsilon / \epsilon^2$ is the value of the optimization problem (\ref{max-min_pb}) at the optimal point.

Let us discuss heuristically why we need to solve this problem, before giving the solution. Note that 
under the null hypothesis our statistic becomes $t_{ij,w} = \sum_{k \in \mathbb{Z}} w_k(\eta^2_{ij,k}-1)$
and it is centered and reduced (due to the normalization $\sum_{k \in \mathbb{Z}} w_k^2 =1/2$).
Under the alternative, 
\begin{eqnarray}
\E_{\bt_{ij}} (t_{ij,w}) = \sum_{k \in \mathbb{Z}} w_k \left( \frac{\theta_{ij,k} }{\ep \sigma_k} \right)^2. \label{bias-alt}
\end{eqnarray}
In order to distinguish the alternative from the null at best, we need to consider the worst parameter ${\bt}_{ij}$ under the alternative and then maximize over possible weights $w_k\geq 0$ verifying the normalization constraints $\sum_k w_k^2=1/2$.

\begin{proposition} \label{sol-opt-contr}
Let $\{\sigma_k\}_{k \in \mathbb{Z}}$ be a sequence of positive real numbers such that $\sigma_k \sim \vert k \vert^s$ as $\vert k \vert$ large enough, for a given $s>0$. Then,    
the optimization problem (\ref{max-min_pb}) has the following solution: 
\begin{eqnarray*}
(\theta_k^\star)^2 & = & v\sigma_k^4 \sqrt{2} \left( 1 - \left( \frac{|k|}{T}\right)^{2\tau}\right)_+, \, T \sim \left(\frac {\kappa_1}{\kappa_2} \right)^{\frac 1{2\tau}} r_\epsilon^{-\frac 1\tau},\,  v = \frac{1}{\kappa_1} \left(\frac{\kappa_2}{\kappa_1}\right)^{\frac{4s+1}{2\tau }} r_\ep ^{2 + \frac{4s+1}\tau} ;\\
w_k^\star & = & \frac{(\theta_k^\star)^2}{2 \sigma_k^2 V_\epsilon}, \mbox{ is such that } \max_{k} w_k^\star \leq r_\ep^{1/(2 \tau)}\to 0;\\
V_\ep &\sim & c(\tau,s) r_\ep^{2 +\frac{4s+1}{2\tau} } , \quad
a(r_\ep) \sim c(\tau,s) \ep^{-2} r_\ep^{2+ \frac{4s+1}{2\tau}},
\end{eqnarray*} 
where the asymptotics are taken as $k \to \infty$ and as $ r_\ep \rightarrow 0$, with
\begin{eqnarray*}
&&c(\tau,s)^2 = 2 (\frac{\kappa_1}{\kappa_2})^{- (4s+1)/ (2\tau)} \frac{ \kappa_3}{\kappa_1^2}, \, 
  \kappa_1= \frac{4 \sqrt{2} \tau}{(4s+1) (4s+2 \tau +1)},\\
&& \kappa_2= \frac{4 \sqrt{2} \tau (2 \pi)^{2 \tau}}{(4s+2 \tau +1) (4s+4 \tau +1)} \mbox{ and } \kappa_3= \frac{1}{4s+1} -\frac{2}{4s+2 \tau+1}  +\frac{1}{4s+4\tau+1},
\end{eqnarray*} 
and where $(x)_+=\max(0,x)$.  
\end{proposition}
The proof of Proposition \ref{sol-opt-contr} is postponed to Appendix. 
Note that $\{w_k^\star\}_k$ and $\{\theta_k^\star\}_k$ check the constraints in (\ref{max-min_pb}), that is,
 $\sum_k (w_k^\star)^2=\frac 1 2$, $
\sum_k (\theta_k^\star)^2 = r_\ep^2(1+o(1))$  and $\sum_k (2 \pi k)^{2 \tau}(\theta_k^\star)^2 =1+o(1), \mbox{ as } r_\ep \to 0
$. 
It is worthwhile to note that due to Proposition \ref{sol-opt-contr} and relation (\ref{bias-alt}), we have
\begin{eqnarray}
\frac 1 2 \sum_k \frac{(\theta_k^\star)^4}{\ep^4 \sigma_k^4 } & = & a^2(r_\ep) \label{cond:a2}  \\
\inf_{\boldsymbol{\theta}_{ij} \in \Sigma (\tau,r_\ep)}   \E_{\bt_{ij}} (t_{ij,w^{\star}}) & =&  a(r_\ep) \label{biais-a}
\end{eqnarray}
and note also that the sequences $\{w_k^\star\}_k$ and $\{\theta_k^\star\}_k$ have a finite number $T$ of non null elements, but $T$ grows to infinity as $r_\ep \to 0$.

\bigskip
 
% Define $t^{\chi^2}$   %the normalized empirical mean of the  $t_{ij}$'s and 
% and $t^{scan}$ as follows,
% \begin{eqnarray}
% t^{\chi^2}& =&  \frac{1}{\sqrt{ MN }} \; \sum_{(i,j) \in I \times J} t_{ij, w^\star}  \label{stat_moderate} \\
% t^{scan}&=&  \max_{\{A \subset I,\, \# A=m, \, B \subset J, \,\# B=n \}}   \frac{1}{\sqrt{ m n}}\; \sum_{(i,j) \in A \times B} t_{ij,w^\star}  \label{stat_max}
% \end{eqnarray}
% where $ w^\star=\{ w^\star_k \}_{k \in \BBz}$ is the sequence of weights which solves the optimization problem (\ref{max-min_pb}). 
%where $t_j$ are given in (\ref{tj}).
Define the test procedures, 
\begin{eqnarray}
\psi^{\chi^2} & =& \1( t^{\chi^2} > H),\mbox{ {\rm with} } t^{\chi^2} =  \frac{1}{\sqrt{ MN }} \; \sum_{(i,j) \in I \times J} t_{ij, w^\star} \label{Test_moderate} \\
\psi^{scan}& =& \1 ( t^{scan} > K ), \mbox{ {\rm with} } \, t^{scan}=  \max_{\xi \in T_{M,N}(m,n)}   \frac{1}{\sqrt{ m n}}\; \sum_{(i,j) \in A_\xi \times B_\xi} t_{ij,w^\star}  \label{Test_high}
\end{eqnarray}
where 
$H$ and $K$ are  positive and $ w^\star=\{ w^\star_k \}_{k \in \mathbb{Z}}$ is the sequence of weights which solves the optimization problem (\ref{max-min_pb}). 
  %Note also that under $\P_0$,
  %the variance of $t_j$ is one since $\sum_{k \in \BBz} w_k^2=\frac 1 2$.

The following theorem gives the upper bounds for the testing rates of the previously defined procedures.

\begin{theorem} \label{thm:upper.bounds} 
Assume (\ref{Sp-C}). 
Suppose that $r_\ep \rightarrow 0$ and recall that 
$$
a(r_\epsilon) \sim c(\tau,s) \ep^{-2} r_\ep^{2+(4s+1)/(2\tau)}.
$$ 
\begin{enumerate}
\item The linear test statistics $\psi^{\chi^{2}}$ defined by (\ref{Test_moderate}) has the following properties.
\begin{description}%Suppose that conditions of Proposition \ref{GD} hold.
\item Type I error probability:  if $H \to \infty$, then $\omega ( \psi^{\chi^2}) = o(1)$.
\item Type II error probability: if  
\begin{equation} \label{cond:chi}
a^2(r_\epsilon)  m n p q \rightarrow  +\infty,
\end{equation}  
choose $H$ such that $H \leq c \cdot a(r_\epsilon) \sqrt{ m n p q}$, for some $0<c<1$,  then $\beta(\psi^{\chi^2},\Theta_{M,N}(\tau,r_\epsilon,m,n))=o(1)$.
 
\end{description}
\item The scan test statistic $\psi^{scan}$ defined by (\ref{Test_high}) has the following properties. \\
Take $K^2 = {2 (1+\delta) ( m\cdot \log(p^{-1})+n\cdot \log(q^{-1}))} $, for some small $\delta >0$, and suppose moreover that $K^2 r_\ep^{1/\tau} /(mn) $ tends to 0 asymptotically. 

%%%%%%%%%%%%PARTIE   ENLEVER APRES VALIDATION%%%%%%%%%%%%%%%%%%%%%%%%%%%
% \bigskip
% (POUR GG, IL FAUT AJOUTER UNE CONDITION : du type $( \log(p^{-1})/n+ \log(q^{-1})/q) =
% o(\epsilon^{-2/(2s+2 \tau +1)})$ pour que la condition $ K^{2}/(mn) (\max_k w_k^{\star})^{2} = o(1)$ soit satisfaite \\
% cette condition implique aussi que $a(r_\epsilon^{\star}) \max_k w_k^{\star} = o(1)$ avec le $r_\epsilon^{\star}$ qui donne 
% le $a(r_\epsilon^{\star})$ sur la 
% boundary detection-- cette derniere condition est necessaire pour calculer la fonction generatrice des moments sous $H_1$) 
% \bigskip
%%%%%%%%%%%%PARTIE   ENLEVER APRES VALIDATION%%%%%%%%%%%%%%%%%%%%%%%%%%%

\begin{description}
\item Type I error:  we have $\omega ( \psi^{scan}) = o(1)$.
\item Type II error: if  
\begin{equation}\label{cond:scan}
\lim \inf \frac{a^2(r_\epsilon) {m n}}{2( m\cdot \log(p^{-1})+n\cdot \log(q^{-1}))} > 1,
\end{equation}  
then $\beta(\psi^{scan},\Theta_{M,N}(\tau,r_\epsilon,m,n))=o(1)$.
\end{description}
\end{enumerate}
%\end{itemize}
\end{theorem}

Consider $\psi$ the test procedure
which combines $\psi^{\chi^2}$ and $\psi^{scan}$ as follows
\begin{eqnarray*}\psi=\max (\psi^{\chi^2},\psi^{scan}).\end{eqnarray*}
 As a consequence of Theorem \ref{thm:upper.bounds}, the test procedure
$\psi$ with $H$ and $K$ properly chosen is such that
$\gamma(\psi,\Theta_{M,N}(\tau,r_\epsilon,m,n))=o(1)$ as soon as either (\ref{cond:chi}) or
(\ref{cond:scan}) hold.

The procedure is rather simple to implement. However, there are difficulties for implementing the
scan procedure. Indeed, computing the scan statistic $t^{scan}$ implies computing standardized sums over
all submatrices of size $m \times n$ in the large matrix $M\times N$. This is computationally infeasible for large values of $M,\, N, m$ and $n$. However, a heuristic algorithm can be implemented as in \cite{BI}, following \cite{SN} and \cite{Shabalinetal}, which is a random procedure finding local maxima. With a sufficiently large choice of random initial values in the algorithm there is a large probability that the algorithm actually finds the global maximum that we aim at.

%%%%%%%%%%%%%%%%%%%%%%%%%%%%%%%%%%%%%%%%%%%%%%%%%%%%%%%%%%
\section{Optimality of the detection boundaries} \label{sec:lower.bounds}
%%%%%%%%%%%%%%%%%%%%%%%%%%%%%%%%%%%%%%%%%%%%%%%%%%%%%%%%%%

We prove here optimality results for the rates that the previous test procedure $\psi$ attained.
However, the optimality is attained under additional hypothesis requiring an 'almost' squared matrix
in the sense that the relative sizes of the submatrix should of the same order in both directions 
(rows and columns sizes). 

\begin{theorem}\label{thm:lower.bound}
Assume (\ref{Sp-C}) and that 
\begin{eqnarray}
\frac{\log \log (p^{-1})}{\log(q^{-1})} \to 0, \quad 
\frac{\log \log (q^{-1})}{\log(p^{-1})} \to 0.
\label{lb1}
\end{eqnarray}
Assume, moreover, that
\begin{eqnarray}
m \cdot \log(p^{-1}) \asymp n \cdot \log(q^{-1}) \label{lb2}
\end{eqnarray}
and that
\begin{eqnarray}  \label{Cond-a-m-n}
\frac{m \cdot \log(p^{-1}) + n \cdot \log(q^{-1})}{mn}=o(1) \ep^{-\frac 2{2\tau +2s+1}}.
\end{eqnarray}

If $r_\ep$ is such that the following conditions are satisfied
\begin{eqnarray}\label{lobo1}
a^2(r_\ep) \cdot m n p q \to 0,
\end{eqnarray}
 \begin{eqnarray}\label{lobo2}
\limsup \frac{a^2(r_\ep)\cdot m n}{ 2(m \cdot \log(p^{-1}) + n \cdot \log(q^{-1}))} < 1,
\end{eqnarray}
then $\inf_\psi \gamma(\psi, \Theta_{M,N}(\tau, r_\ep, m,n)) \to 1$.
\end{theorem}

The proof of Theorem \ref{thm:lower.bound} is given in Section~\ref{sec:lower.bounds}. It follows closely the proof in \cite{BI} with important differences due to the non gaussian likelihoods in this setup.
 
 %%%%%%%%%%%%%%%%%%%%%%%%%%%%%%%%%%%%%
 \subsection{Optimal detection boundary}
 %%%%%%%%%%%%%%%%%%%%%%%%%%%%%%%%%%%%%
 
Theorems \ref{thm:upper.bounds} and \ref{thm:lower.bound} together say that, under assumptions (\ref{Sp-C}),
(\ref{lb1}), (\ref{lb2}) and (\ref{Cond-a-m-n}) a  detection boundary $\tilde{r_\ep}$ is defined via the relations
$$
a^2(\tilde{r_\ep}) \cdot m n p q \asymp 1, \quad a^2(\tilde{r_\ep}) \cdot m n \sim 2 (m \cdot \log(p^{-1})+ n \cdot \log(q^{-1})) .
$$
Therefore, the detection boundary can be written
\begin{equation}\label{sep} \nonumber
a(\tilde{r_\ep}) \asymp \min \left\{ \frac 1{\sqrt{mnpq}}, \, \sqrt{\frac{2( m \cdot \log(p^{-1})+ n \cdot \log(q^{-1}) )}{mn}}\right\},
\end{equation}
with the constant equal to $1$ if $ \displaystyle{\frac{2( m \cdot \log(p^{-1})+ n \cdot \log(q^{-1}) )}{mn} \leq \frac 1{{mnpq}}}$.
By Proposition~\ref{sol-opt-contr}, we have that $a(\tilde{r_\ep}) \sim c(\tau, s) \ep^{-2} (\tilde{r_\ep})^{2 + (4s+1)/(2\tau) }$ as $\ep \to 0$ and $\tilde{r_\ep} \to 0$. It implies that $\ep^{2} a(\tilde{r_\ep}) \to 0$, giving furthermore that one of the following conditions hold:
$$
\ep^{-2} \sqrt{m n p q} \to \infty \mbox{ or } \sqrt{2 ( m \cdot \log(p^{-1})+ n \cdot \log(q^{-1}) )} = o(\ep^{-2} \sqrt{m n}).
$$
The minimax error of the scan test tends to 0, if $K^2 (\tilde{r_\ep})^{1/\tau}/(mn) \to 0$ and a sufficient condition for that is
condition (\ref{Cond-a-m-n}).
%
%$$
%\frac K{\sqrt{mn}} = o\left(\ep^{-\frac 1{2\tau + 2s + 1}} \right).
%$$
%Note that these are not additional assumptions, they merely translate the assumption $r_\ep \to 0$ in Theorem~\ref{thm:upper.bounds} for the sharp detection boundary $\tilde{r_\ep}$. This also explains condition \nref{Cond-a-m-n} in Theorem~\ref{thm:lower.bound}.

We can recover from these results the rates for one dimensional sequences (i.e. $M=N=m=n=1$). In this case it is required that $a(\tilde{r_\ep})$ is asymptotically constant, which means $\tilde{r_\ep} \sim \ep^{4\tau/(4\tau +4s+1)}$ and that is the minimax rate for testing one dimensional signal with Sobolev smoothness $\tau$.

%%%%%%%%%%%%%%%%%%%%%%%%%%%%%%%%%%%%%%%%%%%%%%%%%
%\subsection{Universal statistical procedure}
%%%%%%%%%%%%%%%%%%%%%%%%%%%%%%%%%%%%%%%%%%%%%%%

% \begin{remark}{\bf je me demande si a marcherai comme a: un seul choix de $\tilde w_k^\star$}
% It is clear that our test procedure introduced in the previous part  depends on $r_\epsilon$.  
% One may be interested in dealing with an universal test procedure, i.e., the same statistical procedure 
%  for a range of 
% $r_\epsilon$ and for which it is possible to derive the upper bounds as in Theorem \ref{thm:upper.bounds}.
% In this section, we describe such  a procedure. 

% Denote by   $\tilde{ r}_{ \epsilon}$  the  rates 
% obtained from   ${ a}(\tilde r_{ \epsilon}) $  in \nref{sep}
% using Proposition \ref{sol-opt-contr}   We then 
% consider the test procedures ${\tilde \psi}^{\chi^{2}}$ and   ${\tilde \psi}^{scan}$ given by (\ref{Test_moderate}) and (\ref{Test_high})  but 
% depending on the sequence ${\tilde w}_{k}^{\star}$, solution of the extremal problem (\ref{max-min_pb}) with $\tilde{ r}_{ \epsilon}$ in place of 
% $r_{\epsilon}$.  Put $\tilde \psi = \max\{{\tilde \psi}^{\chi^{2}},\,{\tilde \psi}^{scan}\}$.

% Due to  Proposition 4.1. in \cite{GI}, one gets that any $\boldsymbol{\theta}_{ij} \in \Sigma (\tau, B \tilde{ r}_{ \epsilon})$ with $B \geq 1$,
% \begin{eqnarray}
% \E_{\boldsymbol{\theta}_{ij} } (t_{ij,{\tilde{w}}^{\star}}) &\geq & B^{2} { a}(\tilde r_{ \epsilon}) . \label{univ-biais}
% \end{eqnarray}
% We can prove Theorem~\ref{thm:upper.bounds} for $\tilde \psi$.
% \end{remark}

\subsection{Universal test procedure}

It is clear that our test procedure introduced in the previous part  depends on $r_\epsilon$.  
One may be interested in dealing with an universal test procedure, i.e., a unique  statistical procedure 
 for a range of 
$r_\epsilon$ and for which it is possible to derive the upper bounds as in Theorem \ref{thm:upper.bounds}.
In this section, we describe such  a procedure. 

% Set $${\tilde a}(r_{ \epsilon}) =\left\{ \begin{array}{lcl} (m n p q)^{-1/2} & \mbox{{\rm for }} & \psi^{\chi^2} \\
% \sqrt{2 \frac{m \log (p^{-1})+ n \log (q^{-1})}{mn}} & \mbox{{\rm for }} & \psi^{scan} \end{array}
% \right. $$

Recall that  ${ r}_{ \epsilon}$  can be obtained from ${a}(r_{ \epsilon}) $  
using Proposition \ref{sol-opt-contr}. This gives $r_\ep \asymp (\ep^4 a^2(r_\ep))^{ \tau/(4\tau +4s+1)}$.  Moreover, the optimisation problem (\ref{max-min_pb}) gives associated optimal weights $w_k^\star = w_k^\star(r_\ep)$.
Let us consider $\tilde r_\ep^{\chi^2}$ obtained from $a^2(\tilde r_\ep^{\chi^2}) \asymp (mnpq)^{-1}$ and the associated test procedure ${\tilde \psi}^{\chi^{2}}$ defined in (\ref{Test_moderate}) for the weights $w_k^\star(\tilde r_\ep^{\chi^2}) $. Similarly,  put $\tilde r_\ep^{scan}$ obtained from $a^2(\tilde r_\ep^{scan}) \sim 2(m\log(p^{-1})+n\log(q^{-1}))/(mn)$ and  ${\tilde \psi}^{scan}$ as in   (\ref{Test_high}) with the weights $w_k^\star(\tilde r_\ep^{scan})$.  

The type II error of $\psi^{\chi^2}$ and $\psi^{scan}$ are stated  for 
$\lim a(r_\epsilon)/ {\tilde a}(r_{ \epsilon}) \rightarrow +\infty$ and 
$\lim \inf a(r_\epsilon)/{\tilde a}(r_{ \epsilon}) >1$, respectively, which are equivalent due to Proposition \ref{sol-opt-contr}, 
to   $\lim r_\epsilon/\tilde{ r}_{ \epsilon} \rightarrow +\infty$ and $\lim \inf r_\epsilon/\tilde{ r}_{ \epsilon} >1$. 
Due to  Proposition 4.1. in \cite{GI}, one gets that any $\boldsymbol{\theta}_{ij} \in \Sigma (\tau, B \tilde{ r}_{ \epsilon})$ with $B \geq 1$,
\begin{eqnarray}
\E_{\boldsymbol{\theta}_{ij} } (t_{ij,{w}^{\star}(\tilde{ r}_{ \epsilon})}) &\geq & B^{2} a(\tilde{ r}_{ \epsilon}). \label{univ-biais}
\end{eqnarray}

Taking $H$  and $K$ as  in Theorem \ref{thm:upper.bounds} and  using  (\ref{univ-biais}),
we derive  for ${\tilde \psi}^{\chi^{2}}$ and ${\tilde \psi}^{scan}$ the same upper bounds as in Theorem \ref{thm:upper.bounds}. 
For ${\tilde \psi}^{scan}$, note that 
the condition $(K/(\sqrt{ nm} ) \max w_k^{\star} (\tilde r_\ep^{scan})= o(1)$ is satisfied as soon as 
$(m \log(p^{-1}) + n \log(q^{-1}))/(nm) =o(\epsilon^{-2/(2s+1+2\tau)})$.

%%%%%%%%%%%%%%%%%%%%%%%%%%%%%%%%%%%%%
\section{Proofs} \label{sec:Proofs}
%%%%%%%%%%%%%%%%%%%%%%%%%%%%%%%%%%%%%
Let us start with a preliminary result that gives an approximation of the  moments generating function of $t_{ij,w^{\star}}$ under $H_0$. 
%%% PAS BESOIN DE LA GENERATRICE SOUS L'ALTERNATIVE %%%%%%%%%%%%%%%%%%%%%%%%%%%%%%%%%%%%%%%%%%
\begin{lemma} \label{generatrice-function} 
For any real number $\lambda$ such that $\lambda \max_k w_k^{\star}=o(1)$, 
\begin{eqnarray*}
\mathbb{E}_0 (\exp (\lambda t_{ij,w^{\star}})) & =& \exp \left( \frac{\lambda^{2}}{2} (1+o(1)) \right) \label{generatrice-H0}.
\end{eqnarray*}
%
%assume moreover that $r_\epsilon \rightarrow 0$ and $a(r_\epsilon)\max_k w_k^{\star}=o(1)$,  
%%(FAUT VOIR CETTE CONDITION ET SON LIEN AVEC DISCUSSION APRES THM 3.1)  
%then
%\begin{eqnarray}
%\mathbb{E}_{\boldsymbol{\theta}_{ij}} (\exp (\lambda t_{ij,w^{\star}})) & =& \exp\left (  \lambda a(r_\epsilon)(1 + o(1))+ \frac{\lambda^{2}}{2} (1+o(1)). \right) \label{generatrice-H1}
%\end{eqnarray} 
\end{lemma}
The proof of Lemma \ref{generatrice-function} is postponed in the appendix.

%%%%%%%%%%%%%%%%%%%%%%%%%%%%%%%%%%%%%%%%%%%%%%%%%%%%%%%
\subsection{Proof of Theorem~\ref{thm:upper.bounds}}
%%%%%%%%%%%%%%%%%%%%%%%%%%%%%%%%%%%%%%%%%%%%%%%%%%%%%%%

%\subsection{Preliminaries}
Observe that under $H_0$, $t_{ij,w}$ are i.i.d. random variables with  zero mean and unit variance. 
Indeed, %setting ${\tilde \eta}_{ij,k}=\displaystyle{ \frac{\eta_{ij,k}}{\sigma_{k}}}$,  
one gets $\V_0 (t_{ij,w}) = \sum_k w_k^2 \V({\eta}_{ij,k}^2)/\sigma_k^2 = 2 \sum_k w_k^2 = 1$.
Under the alternative, for all $\boldsymbol{\theta}_{ij}  \in  \Sigma (\tau,r_{\epsilon})$, 
\begin{eqnarray*}
\E_{\boldsymbol{\theta}_{ij}}(t_{ij,w}) 
& =&  \sum_k w_k  \frac{\theta_{ij,k}^{2}}{\sigma^2_{k}\epsilon^{2}} \nonumber \\
%& =& \xi_{ij} \kappa(\boldsymbol{\theta}_{ij}, w) \frac{1}{\epsilon^{2}} \nonumber \\
\V_{\boldsymbol{\theta}_{ij}}(t_{ij,w})
& =& \sum_k w_k^{2}  (2 + 4  \frac{\theta_{ij,k}^{2}}{\sigma_{k}^{2}\epsilon^{2}} ) \nonumber \\ %\xi_{ij}
%& =&   1 +4  \sum_k w_k^{2} \frac{\theta_{ij,k}^{2}}{\sigma_{k}\epsilon^{2}} \nonumber \\% \xi_{ij}
& \leq &  1 +4  \sup_k w_k \cdot \sum_k w_k \frac{\theta_{ij,k}^{2}}{\sigma_{k}^{2}\epsilon^{2}}  %\nonumber \\% \xi_{ij}
%& =& 
=1 +4  \sup_k w_k \cdot \E_{\boldsymbol{\theta}_{ij}}(t_{ij,w}) . \nonumber
\end{eqnarray*}

Due to the previous relations, for any $\overline{\boldsymbol{\theta}} \in \Theta_{M,N}(\tau,r_\epsilon,m,n)$
\begin{eqnarray}
\E_{\overline{\boldsymbol{\theta}}} (t^{\chi^2})& =& \frac{1}{\sqrt{ MN}} \; \sum_{(i,j) \in A_\xi \times B_\xi} \E_{\boldsymbol{\theta}_{ij}} (t_{ij,w^\star}) \nonumber \\
%& =& \frac{1}{\sqrt{ MN}} \; \sum_{(i,j) \in I \times J} \sum_k w_k  \xi_{ij}\frac{\theta_{ij,k}^{2}}{\epsilon^{2}} \nonumber \\
& \geq & \sqrt{M N} p q \cdot a(r_\epsilon) =  \sqrt{m n p q} \cdot a(r_\epsilon), \label{mean-Chi2-alt} 
\end{eqnarray}
where  the penultimate inequality follows from (\ref{biais-a}).
Moreover, for the variance we have
\begin{eqnarray}
\V_{\overline{\boldsymbol{\theta}}} (t^{\chi^2})
& =& \frac{1}{ MN} \; \sum_{(i,j) \in I \times J} \V_{\overline{\boldsymbol{\theta}}_{ij}} [t_{ij,w^\star} (\1(\xi_{ij}=1)+\1(\xi_{ij} \not =1))] \nonumber \\
&\leq & 1+ \frac 4{MN} \sum_{(i,j) \in A_\xi \times B_\xi}\sum_k (w_k^\star)^2 \frac{\theta_{ij,k}^2}{\ep^2 \sigma_k^2}
\nonumber\\
%& = &\frac{1}{ MN} \; \sum_{(i,j) \in I \times J}\left( 1 + 4 \sum_k w_k^{2}  \xi_{ij}\frac{\theta_{ij,k}^{2}}{\epsilon^{2}} \right) \nonumber \\
& \leq & 1 + 4 \max_k w_k^\star  \frac{1}{ \sqrt{MN}} \E_{\overline{\boldsymbol{\theta}}} (t^{\chi^2}). \label{var-Chi2-alt} 
\end{eqnarray}
Recall that  $\displaystyle{\max_k w_k^\star \stackrel{r_\epsilon \rightarrow 0}{\longrightarrow} 0}$ (see Proposition~\ref{sol-opt-contr}).  

%
%Similarly, one obtains, 
%\begin{eqnarray}
%\E_{\overline{\boldsymbol{\theta}}} (t^{scan})& \geq & \frac{1}{\sqrt{ MN}} \; \sum_{(i,j) \in I \times J} \E_{\boldsymbol{\theta}_{n,m}} t_{ij} \nonumber \\
%& =& \frac{1}{\sqrt{ MN}} \; \sum_{(i,j) \in I \times J} \sum_k w_k  \xi_{ij}\frac{\theta_{ij,k}^{2}}{\epsilon^{2}} \nonumber \\
%& \geq & \sqrt{N M} p q a(r_\epsilon). \label{mean-Chi2-alt} 
%\end{eqnarray}
%
%
%\begin{eqnarray}
%\V_{\overline{\boldsymbol{\theta}}} (t^{\chi^2})& =& \frac{1}{ MN} \; \sum_{(i,j) \in I \times J} \V_{\overline{\boldsymbol{\theta}}_{ij}} t_{ij} \nonumber \\
%& = &\frac{1}{ MN} \; \sum_{(i,j) \in I \times J}\left( 1 + 4 \sum_k w_k^{2}  \xi_{ij}\frac{\theta_{ij,k}^{2}}{\epsilon^{2}} \right) \nonumber \\
%& \leq & 1 + 4 \sup_k w_k  \frac{1}{ \sqrt{MN}} \E_{\overline{\boldsymbol{\theta}}} (t^{\chi^2}) \label{var-Chi2-alt} 
%\end{eqnarray}

{\bf Type I error probability of $\psi^{\chi^2}$.}
Since $\max_k w_k^{\star}=o(1)$, the asymptotic standard normality of $t^{\chi^2}$ under the null  follows  from Lemma 3.1 in \cite{IS.02a}
then,
\begin{eqnarray*}
\P_0( t^{\chi^2} >H) & =& \Phi(-H)+ o(1),
\end{eqnarray*}
where $\Phi$ stands for the cdf of a standard Gaussian random variable.

 \vspace{0.1cm} %\noindent
 
{\bf Type II error probability of $\psi^{\chi^2}$ uniformly over  $\Theta_{M,N}(\tau,r_\epsilon,m,n)$.} % for all $\theta$ in $\Theta_d(\tau,r_\epsilon,b)$.
We deduce from (\ref{var-Chi2-alt}) that  
 $\Var_{\overline{\boldsymbol{\theta}}} (t^{\chi^2})= 1 + o (\E_{\overline{\boldsymbol{\theta}}} (t^{{\chi^2}}))$, uniformly over  ${\overline{\boldsymbol{\theta}}}\in \Theta_{M,N}(\tau,r_\epsilon,n,m)$. 
 % provided that $\E_{\overline{\boldsymbol{\theta}}} (t)  \rightarrow +\infty$.

For all  $\overline{\boldsymbol{\theta}}$ in $\Theta_{M,N}(\tau,r_\epsilon,m,n)$, by Markov's inequality and relation (\ref{mean-Chi2-alt}),
\begin{eqnarray*}
\P_{\overline{\boldsymbol{\theta}}} ( t^{\chi^2} \leq H) 
& \leq &
\P_{\overline{\boldsymbol{\theta}}}  \left( \vert t^{\chi^2}  - \E_{\overline{\boldsymbol{\theta}}} (t^{\chi^2})\vert \geq   \E_{\overline{\boldsymbol{\theta}}} (t^{\chi^2}) - H \right)\nonumber \\
& \leq & \frac{ \Var_{\overline{\boldsymbol{\theta}}} (t^{\chi^2})}{(\E_{\overline{\boldsymbol{\theta}}} (t^{\chi^2})- H )^2}\\
& \leq & \frac{1 + 4 \max_k w_k^{\star} \E_{\overline{\boldsymbol{\theta}}} (t^{\chi^2})/\sqrt{MN}}{(\E_{\overline{\boldsymbol{\theta}}} (t^{\chi^2}) -H)^2}\\
&\leq & \frac 1{(1-c)^2 \E_{\overline{\boldsymbol{\theta}}} (t^{\chi^2})^2}
+ \frac {4 \max_k w_k^\star}{(1-c)^2 \sqrt{MN} \E_{\overline{\boldsymbol{\theta}}} (t^{\chi^2})}
=o(1),
\end{eqnarray*}
provided that $ a(r_\epsilon)\sqrt{ m n p q }\rightarrow +\infty$ and $H \leq c\cdot a(r_\epsilon)\sqrt{m n p q}$ for some $0<c<1$.

{\bf Type I error probability of $\psi^{scan}$.}

Under the assumption~(\ref{Sp-C}), we can check that
\begin{eqnarray*}
\log(C_M^m \cdot C_N^n) \sim m \cdot \log(p^{-1}) + n\cdot \log(q^{-1}). \label{equi-Comb}
\end{eqnarray*}
Applying Markov's inequality,
\begin{eqnarray}
 \P_0(t^{scan} > K) & \leq &\sum_{\xi \in T_{M,N}(m,n)}    
 \P_0( \frac{1}{\sqrt{ m n}}\; \sum_{(i,j) \in A_\xi \times B_\xi} t_{ij,w^{\star}}  >K) \nonumber \\
& = &C^{m}_{M} C^{n}_{N}  \P_0( \frac{1}{\sqrt{ m n}}\; \sum_{(i,j) \in A_\xi \times B_\xi} t_{ij,w^{\star}}  >K), \nonumber \\
& \leq & C^{m}_{M} C^{n}_{N}  e^{-K^2} \E_0 (\exp(\sum_{(i,j) \in A_\xi \times B_\xi} \frac{K}{\sqrt{m n}} t_{ij,w^{\star}})) \nonumber\\
& \leq & C^{m}_{M} C^{n}_{N}  e^{-K^2} \left(\E_0(\exp(\frac{K}{\sqrt{m n}} t_{11,w^{\star}}))\right)^{m n}.\label{interm}
\end{eqnarray}
Set $\lambda={K}/{\sqrt{m n}}$ with $K = \sqrt{2(1+\delta) \log(C_N^n C_M^m)}$, for some small $\delta>0$  and note that $K/\sqrt{mn} \max_k w_k^{\star} \leq K r_\ep^{1/(2\tau)} /\sqrt{mn} =o(1)$ by assumption in our theorem; 
then, applying   Lemma \ref{generatrice-function} we obtain that 
%
%Recall that ${\tilde \eta}_{ij,k}= \frac{\eta_{ij,k}}{\sigma_{ij,k}}$ is a standard Gaussian real random variable. 
% As $t_{11,w^{\star}}=\sum_k w_k^{\star} ({\tilde \eta}_{11,k}^2 -1)$ is a sum of independent $\chi^2$ random variables, we use the fact that
$$\mathbb{E}_0(\exp(\lambda t_{11,w^{\star}}))  =\exp(\frac{\lambda^2}{2}(1+o(1))).$$
 %$\max_k w_k^{\star}=o(1)$ as $r_\ep \rightarrow 0$ (see Lemma~\ref{sol-opt-contr}).
Next,  by plugging $(\E_0(\exp(\frac{K}{\sqrt{nm}} t_{11,w^{\star}})))^{nm} = \exp( \frac{K^2}{2nm} (1+o(1)) )$
 into (\ref{interm}), we obtain
\begin{eqnarray*}
\P_0(t^{scan} > K) 
&\leq &
C^{m}_{M} C^{n}_{N}  e^{-K^2/2 (1+o(1))} =o(1),
\end{eqnarray*}
due to the choice of $K = \sqrt{2(1+\delta) \log(C_N^n C_M^m)}$, for some small $\delta>0$.

% To end this part, let us note that  condition $\lambda \;\max_k w_k^{\star}=o(1)$ is satisfied since  
% $\max_k w_k^{\star}=o(1)$ as $r_\ep \rightarrow 0$ (see Lemma~\ref{sol-opt-contr}) and  as $n$, $m$, $M$ and $N$ go to infinity, using Stirling's approximation, we obtain 
% \begin{eqnarray*}
% \frac{K^{2}}{nm} & \asymp &  \frac{\log M }{n} + \frac{\log N}{m}.
% \end{eqnarray*}
% The right-hand side  of this latter relation goes to zero if 
% $N \asymp M$ and if $n \asymp N^{1-b_1}$ and $n \asymp N^{1-b_2}$ where $b_1 \in (0,1)$ and $b_2 \in (0,1)$. 
 
%t^{scan}&=&  \max_{\{A \subset I, B \subset J \}}   \frac{1}{\sqrt{ m n}}\; \sum_{(i,j) \in A \times B} t_{ij} 

 \vspace{0.1cm} %\noindent
 
{\bf Type II error probability of $\psi^{scan}$ uniformly over $\Theta_{M,N}(\tau,r_\epsilon,m,n)$.}  For any $\overline{\boldsymbol{\theta}} \in \Theta_{M,N}(\tau,r_\epsilon,m,n)$, it exists $A \subset I$ and $B \subset J$ such that 
$ \# A=m$, $ \# B=n$ and $\xi_{ij}=\1((i,j) \in A\times B)$;   using the inequality $t^{scan} \geq \frac{1}{\sqrt{ m n}}\; \sum_{(i,j) \in A \times B} t_{ij,w^{\star}}$, 
  we obtain
 \begin{eqnarray}
 \P_{\overline{\boldsymbol{\theta}}} (t^{scan} \leq K) %& =& 
%\P_{\overline{\boldsymbol{\theta}}}(\max_{\{A \subset I, B \subset J \}} \frac{1}{\sqrt{ m n}}\; \sum_{(i,j) \in A \times B} t_{ij,w^{\star}} \leq  K)  \nonumber \\
& \leq & \P_{\overline{\boldsymbol{\theta}}} (\frac{1}{\sqrt{ m n}}\; \sum_{(i,j) \in A \times B} t_{ij,w^{\star}} \leq  K)  \nonumber \\
& \leq & \frac{ \V_{\overline{\boldsymbol{\theta}}} (\frac{1}{\sqrt{ m n}}\; \sum_{(i,j) \in A \times B} t_{ij,w^{\star}})}{(\E_{\overline{\boldsymbol{\theta}}} (\frac{1}{\sqrt{ m n}}\; \sum_{(i,j) \in A \times B} t_{ij,w^{\star}} )- K )^2} \nonumber .
\end{eqnarray}
Due to (\ref{biais-a}), we have 
\begin{eqnarray*}
\E_{\overline{\boldsymbol{\theta}}} (\frac{1}{\sqrt{ m n}}\; \sum_{(i,j) \in A \times B} t_{ij,w^{\star}})
& = & \frac{1}{\sqrt{ m n}}\; \sum_{(i,j) \in A \times B} \E_{\boldsymbol{\theta}_{ij}} (t_{ij,w^{\star}}) \geq a(r_\ep) \sqrt{mn}.
\end{eqnarray*}
By assumption \nref{cond:scan} we have $\liminf a(r_\ep) \sqrt{mn} /K \geq (1+\delta)^{-1/2}$, which implies that $K \leq a(r_\ep) \sqrt{mn (1+\delta)}/(1+\tilde \delta)$ for some $\tilde \delta >0$ and then $K \leq c a(r_\ep) \sqrt{mn} \leq c \E_{\overline{\boldsymbol{\theta}}} (\frac{1}{\sqrt{ m n}}\; \sum_{(i,j) \in A \times B} t_{ij,w^{\star}})$ for some $0<c<1$ if $\delta$ is small enough.

Now, acting as for getting Equation (\ref{var-Chi2-alt}), we have
\begin{eqnarray*}
\V_{\overline{\boldsymbol{\theta}}} (\frac{1}{\sqrt{ m n}}\; \sum_{(i,j) \in A \times B} t_{ij,w^{\star}}) &=&
\frac{1}{ m n}\; \sum_{(i,j) \in A \times B} \V_{\boldsymbol{\theta}_{ij}}(t_{ij,w^{\star}})\\
&\leq & 1+\frac{4 \max_k w_k^\star}{mn}\sum_{(i,j) \in A \times B} 
\E_{{\boldsymbol{\theta}}_{ij}} (  t_{ij,w^{\star}} )\\
&\leq & 1+\frac{4 \max_k w_k^\star}{\sqrt{mn}}
\E_{\overline{\boldsymbol{\theta}}} (\frac{1}{\sqrt{ m n}}\; \sum_{(i,j) \in A \times B} t_{ij,w^{\star}} ).
\end{eqnarray*}
Finally,
$$
 \P_{\overline{\boldsymbol{\theta}}} (t^{scan} \leq K)
 \leq \frac 1{(1-c)^2 a^2(r_\ep)mn} + \frac{4\max_k w_k^\star}{(1-c)^2 a(r_\ep) mn} = o(1).
$$

%%%%%%%%%%%%%%%%%%%%%%%%%%%%%%%%%%%%%%%%%%%%%%%%%%%%
\subsection{Proof of Theorem~\ref{thm:lower.bound}} \label{sec:lower.bound}
%%%%%%%%%%%%%%%%%%%%%%%%%%%%%%%%%%%%%%%%%%%%%%%%%%%%%%

%The prior we consider is a  classical one for the functional Gaussian model and it is referenced in section 4.3 of  \cite{IS.02a} as
%the symmetric Three-point Factors. \\

The usual steps for proving the lower bounds are the following. First, we bound from below the minimax total error probability  by reducing the set of parameters. Next, we choose a prior on the reduced set of parameters and bound the testing risk from below with a Bayesian risk. Finally, this Bayesian risk is large if a $\chi^2$-distance between the likelihoods under the null and under the mixture of alternatives is small.

Recall that  $\{\theta_k^{\star}\}_{k \in \mathbb{Z}}$ is the solution of the optimisation problem (\ref{max-min_pb}) and let 
 us choose a matrix $\xi$ in the set $T_{M,N}(m,n)$, $\xi=\I((i,j) \in A \times B)$ where $A=A_\xi$ is a set of $m$ rows out of $M$ and $B=B_\xi$ a set of $n$ columns out of $N$. Denote by 
$$
\mathcal{T}_{M,N}(\tau,r_\ep,m,n) = \{\overline{\boldsymbol{\theta}} = [\xi_{ij} \{\pm \theta_{k}^\star\}_k]_{(i,j)\in I\times J}, \xi \in T_{M,N}(m,n) \}.
$$
This is the reduced set of parameters, i.e., a subset of the alternative $\Theta_{M,N}(\tau,r_\epsilon,m,n)$ in our test.

A prior measure on the reduced set will choose $\xi$ with equal probability in the set $T_{M,N}(m,n)$; given  $\xi$, 
the  $(\theta_{ij})$'s associated with non-zero $\xi_{ij}$ are i.i.d. and  
 for $(i,j)$ such that $\xi_{ij}=1$, the prior  will choose with equal probability between $\theta_k^{\star}$ and $-\theta_k^{\star}$, independently for each $k$. 
We can write $\pi_{ij,k} = \frac 12 (\delta_{-\theta_k^{\star}}+\delta_{\theta_k^{\star}})$, where $\delta$ stands for the Dirac measure,  and $\pi_{ij} = \prod_k \pi_{ij,k} $. Let us define 
$$
\overline \pi = \frac 1{C_M^m C_N^n} \sum_{\xi \in T_{M,N}(m,n)} \prod_{(i,j)\in A_\xi \times B_\xi}\pi_{ij} %I(\xi_{ij}=1),
$$ 
the global prior on $\overline{\boldsymbol{\theta}}$'s in $\mathcal{T}_{M,N}(\tau,r_\ep,m,n)$. 

Let us write the likelihood ratio of one active component, i.e., when $(i,j)$ is such that $\xi_{ij}=1$,
$$
\frac{d \P_{\pi_{i j}}}{d \P_0}(\{x_{i j,k}\}_k) 
= \prod_k \exp(-\frac{{\theta_k^{\star}}^2}{2\epsilon^2 \sigma_k^2}) 
\cosh (x_{i j,k} \frac{\theta_k^{\star}}{\epsilon^2 \sigma_k^2}).%, \mbox{ if } \xi_{i j}=1.
$$

Set $\overline X= [\{ x_{ij,k} \}_k]_{(i,j)}$. 
Then the likelihood ratio with respect to the null hypothesis of our observations becomes:
\begin{equation*}
L_{\overline \pi}(\overline X) = \frac{d\P_{\overline \pi}}{d\P_0}([\{ x_{ij,k} \}_k]_{(i,j)})
=\frac 1{C_M^m C_N^n} \sum_{\xi \in T_{M,N}(m,n)} 
\frac{d\P_{\xi}}{d\P_0}(\overline X), 
\end{equation*}
where
$$
\frac{d\P_{\xi}}{d\P_0}(\overline X)=
\prod_{(i,j) \in A_\xi \times B_\xi} \frac{d\P_{\pi_{ij}}}{d\P_0}(\{x_{ij,k}\}_k).
$$

In order to prove indistinguishability, we see that
\begin{eqnarray*}
\gamma & = & \inf_{\psi \in [0,1]} (w(\psi)+ \sup_{\overline{\boldsymbol{\theta} }\in \Theta_{M,N}(\tau,r_\epsilon,m,n)} \E_{\overline{\boldsymbol{\theta}}} [1-\psi(\overline{X})] ) \\
& \geq & \inf_{\psi \in [0,1]} (w(\psi)+ \sup_{\overline{\boldsymbol{\theta} } \in \mathcal{T}_{M,N}(\tau,r_\epsilon,m,n)} \E_{\overline{\boldsymbol{\theta}}} [1-\psi(\overline{X})] )\\
& \geq & \inf_{\psi \in [0,1]} (w(\psi)+ \sum_{\overline{\boldsymbol{\theta}}  \in \mathcal{T}_{M,N}(\tau,r_\epsilon,m,n)} \overline{\pi} (\overline{\boldsymbol{\theta}}) \E_{\overline{\boldsymbol{\theta}}} [1-\psi(\overline{X})] )\\
& \geq & \inf_{\psi \in [0,1]} (\E_0(\psi(\overline{X}))+ \E_0[(1-\psi(\overline{X})) L_{\overline \pi}(\overline X)]).
\end{eqnarray*}
This infimum is attained for the likelihood ratio test $\psi^\star(\overline X) = \1(L_{\overline \pi}(\overline X) >1)$. By Fatou lemma, we have
$$
\liminf \gamma \geq \E_0(\liminf \psi^\star(\overline{X}) + (1-\psi^\star(\overline{X}))L_{\overline \pi}(\overline X)),
$$
which implies that $\gamma \to 1$ if $L_{\overline \pi}(\overline X) \to 1$ in $\P_0$-probability.
In order to prove this sufficient condition, it is enough to check that
\begin{equation*}\label{cond:indisting}
\E_0(L_{\overline \pi}(\overline X)^2) \leq 1+o(1).
\end{equation*}

However, this last inequality can not be obtained; indeed, too many events with  small probability are summed up in the expected value of the square likelihood ratio. Therefore, we modify slightly the likelihood ratio, by truncation, as follows:
\begin{equation*}
\tilde L_{\overline \pi}(\overline X) 
=\frac 1{C_M^m C_N^n} \sum_{\xi \in T_{M,N}(m,n)} 
\prod_{(i,j) \in A_\xi \times B_\xi} \prod_k \exp(-\frac{{\theta_k^{\star}}^2}{2\epsilon^2 \sigma_k^2}) 
\cosh (x_{i j,k} \frac{\theta_k^{\star}}{\epsilon^2 \sigma_k^2}) \1(\Gamma_\xi),
\end{equation*}
where the event $\Gamma_\xi$ is defined for some small $\delta_1>0$ as follows

\begin{eqnarray}\label{gamma_xi}
\Gamma_{\xi} &= & 
\left\{\frac 1{a(r_\epsilon) \sqrt{hl}} \sum_{(i,j)\in A_V \times B_V} \sum_k 
\left( \log \cosh(u_k \cdot \frac{x_{ij,k}}{\epsilon \sigma_k}) -\frac{u_k^2}2 + \frac{u_k^4}4 \right) 
\leq T_{hl},  \right. \\
&& \left. \forall  V \in T_{M,N}(h,l) 
\mbox{ {\rm such that} } A_V \subset A_\xi, \;  B_V \subset B_\xi \right., \nonumber \\  
&& \left.   \mbox{ {\rm and } }  \;    \forall h, \; \forall l  \mbox{ {\rm such that} } \; \delta_1 m \leq h \leq m \mbox{ and } \delta_1 n \leq l\leq n
\right\} , \nonumber
\end{eqnarray}
with $u_k = \theta_k^\star/(\epsilon \sigma_k)$ and 
$$
T_{hl}^2 = 2 (\log(C_M^h C_N^l)(1+A_\ep) + \log(mn)),
$$ 
for small $A_\ep=o(1)$. Moreover, in order to finish the proof of the lower bounds, Lemma~\ref{lemma:aux.lower.bound},~1.
will require $A_\ep$ to be such that
\begin{equation}\label{Aep}
A_\ep \cdot \log(C_M^h C_N^l) \to \infty \mbox{ and } A_\ep /(a(r_\ep) \max_k w_k^\star) \to \infty.
\end{equation}
We suggest to take, e.g. $A_\ep$ the largest value between $A_{1,\ep} = 1/\log((a(r_\ep) \max_k w_k^\star)^{-1}) $ and $A_{2,\ep} = 1/\log(C_M^h C_N^l)$.

Let us see that, as $r_\ep < \tilde r_\ep$,
\begin{eqnarray*}
a(r_\ep) \cdot \max_k w_k^\star 
&\leq & O(1) a(\tilde r_\ep) \tilde r_\ep^{\frac 1{2\tau}} \\
&\leq & O(1) a(\tilde r_\ep)^{1+ \frac1{4\tau +4s+1}} \ep^{\frac 2{4\tau +4s+1}}\\
&\leq & O(1) (a(\tilde r_\ep) \cdot \ep^{\frac 1{2\tau +2s +1}})^{\frac{4\tau +4s +2}{4\tau +4s +1}} = o(1),
\end{eqnarray*}
by assumption~\nref{Cond-a-m-n}. This also implies that $A_{1,\ep} = o(1)$ and that $A_\ep = o(1)$.

In the following we shall denote by $V \subset \xi$ any matrix of size $M \times N$ 
such that  $V = \1((i,j) \in A_V \times B_V)$, with $A_V \subset A_\xi $ and $B_V \subset B_\xi$. Note that $V$ may have value 1 only in a submatrix where $\xi$ has value 1. 

The idea is that the random variable in this event is truncated at the values predicted by large deviations and this is sufficient to diminish the second-order moment of the likelihood ratio. 

Let us denote $\Gamma = \displaystyle{ \cap _{\xi \in T_{M,N}(m,n)} }\Gamma_\xi$. Then, for some fixed $\delta>0$, let us consider  the event ${\cal E} = \{|L_{\overline \pi}( \overline X) - 1| >\delta \}$
\begin{eqnarray*}
\P_0({\cal E}) & = & \P_0 ({\cal E} \cap \Gamma) + \P_0 ({\cal E} \cap \Gamma ^C)\\
& \leq &  \delta^{-2} \E_0 ((L_{\overline \pi}( \overline X) - 1)^2 \1(\Gamma))
+ \P_0(\Gamma^C)\\
& \leq & \delta^{-2}[(\E_0(\tilde L_{\overline \pi}(\overline X) ^2)-1) -2(\E_0(\tilde L_{\overline \pi}(\overline X)) -1)+(\P_0(\Gamma)-1) ]+ \P_0(\Gamma^C),
\end{eqnarray*}
where $\Gamma^C$ denotes the complement of $\Gamma$. 
Then, it remains to prove the following lemma to complete  the proof of Theorem~\ref{thm:lower.bound}.
\begin{lemma}\label{lemma:aux.lower.bound}
Under the assumptions of Theorem~\ref{thm:lower.bound} we have the following:
\begin{enumerate}
\item $\P_0(\Gamma) \to 1$.

\item $\E_0(\tilde L_{\overline \pi}(\overline X) ) \to 1$.

\item $\E_0(\tilde L_{\overline \pi}(\overline X) ^2) \leq 1+ o(1)$.
\end{enumerate}
\end{lemma}
The proof of Lemma \ref{lemma:aux.lower.bound} is postponed  in Section~\ref{sec:aux.lower.bound}.

\section{Appendix} \label{sec:appendix}

\subsection{Auxiliary results for the lower bounds} \label{subsec:aux.results}

\begin{lemma} \label{lem:stoch-order} As  $N \rightarrow  +\infty$,
 $n \rightarrow  + \infty$ 
 and $n/N \rightarrow 0$, if ${\tilde q}=\displaystyle{ \frac{2n}{N-n}}$ then 
$$
\P_{\mathcal{HG}(N,n,n)}(X\geq k) \leq \P_{\mathcal{B}(n,{\tilde q})}(Y\geq k),  
$$
for all integer numbers $k$, where $\mathcal{HG}(N,n,n)$ and $\mathcal{B}(n,{\tilde q})$ denote the hypergeometric distribution and the Binomial distribution, respectively. \end{lemma} 
\noindent
{\bf Proof of Lemma \ref{lem:stoch-order}} 
The proof is based on Theorem 1 (d) in \cite{KM.10} that states that for positive integers $M$, $A$, $m$, $n$  
such that $m \leq M$ and for $q \in (0,1]$, 
$$
\P_{\mathcal{HG}(M,A,m)}(X\geq k) \leq \P_{\mathcal{B}(n, q)}(Y\geq k)
$$ for all integers $k$ if and only if $\inf (m, A) \leq n$ and 
${C_{M-A}^m}/{C_M^m} \geq (1-q)^{n}$ . 

In our case, it is sufficient to check that asymptotically 
${C_{N-n}^n}/{C_N^n} \geq (1-{\tilde q})^{n}$.
As $N \rightarrow \infty $, $n \rightarrow \infty$ and $n/N \rightarrow 0$, 
using Stirling's formula and after simplifications, one gets, 
\begin{eqnarray*}
\frac{C_{N-n}^n}{C_N^n} &  =&
\frac{((N-n)!)^{2}}{N! (N-2n)!} \nonumber \\
& \sim & \exp \left((2(N-n) +1)\log (N-n) - (N+\frac 12)\log (N) - (N-2n+\frac 12) \log (N-2n)\right) \nonumber \\ 
&\sim & \exp\left( (2(N-n) +1)\log (1- \frac nN) - (N-2n+\frac 12) \log (1-\frac{2n}N)  \right) \nonumber \\
& \sim & \exp\left(-2\frac{n^2}{N} +o(\frac{n^2}N)\right), 
\end{eqnarray*}
which is larger than $(1-{\tilde q})^{n}$ provided that  ${\tilde q} \geq \frac{2n}{N}$;  this is satisfied with ${\tilde q}=\displaystyle{\frac{2n}{N-n}}$.
 \hfill \endproof
%%%%%%%%%%%%%%%%%%%%%%%%%%%%%%%%%%%%%%%%%%%%%%%%%
\begin{lemma}\label{generating_cosh} If $r_\ep \to 0$ such that $a(r_\ep) \cdot \max_k w_k^{\star} = o(1)$, then for any $\lambda>0$ such that  $\lambda =O(1)$,
$$
\E_0 (\exp(\lambda \sum_k Z_{11,k})) = \exp(\frac{\lambda^2 a^2(r_\ep)}{2}(1+O(a(r_\ep) w_k^{\star}))).
$$
\end{lemma}

\noindent
{\bf Proof of  \ref{generating_cosh}}  Let us see that for bounded $\lambda > 0$, for $u\to 0$ and a standard Gaussian random variable $\eta$, we have:
$$
\E(\exp(\lambda \cdot (\log \cosh(u \cdot \eta) - \frac{u^2}2 +\frac{u^4}4 + O(u^6)))) = \exp(\frac{\lambda^2 u^4}{4} + O(u^6)).
$$ 
This proof can be adapted in \cite{GI} (cf Lemma A.1).  Now, we apply this result for each $k$ and 
note that $u_k^{2}= a(r_\epsilon)\cdot w_k^{\star} \leq a(r_\epsilon)\cdot \max_k w_k^{\star}=o(1)$ by assumption. 
Using  $\sum_k u_k^4 = 2 a^2(r_\ep)$, we get
\begin{eqnarray*}
\E_0 (\exp(\lambda \sum_k Z_{11,k})) &=& \exp ( \frac{\lambda^2 }{4} \sum_k u_k^4 (1 + O(u_k^2)))\\
&=& \exp(\frac{\lambda^2 a^2(r_\ep)}{2} (1+O(a(r_\ep) \cdot \max_kw_k^\star))).
\end{eqnarray*}
$\;$ \hfill  \endproof

\subsection{Proof of Lemma~\ref{lemma:aux.lower.bound}}\label{sec:aux.lower.bound}

Take a small $\delta >0$. The detection boundary $a(r_\ep)$ satisfies (\ref{lobo2}), so the most difficult case is when the limit is close to 1. Therefore, we shall assume that
\begin{equation}\label{asim} \nonumber
a^2(r_\ep) mn \sim (2-\delta) (m \cdot \log(p^{-1})+n\cdot \log(q^{-1})).
\end{equation}
This implies
\begin{equation}\label{aasymp}
a^2(r_\ep) \asymp \frac{\log(p^{-1})}n +\frac{\log(q^{-1})}m.
\end{equation}

Let us see that the random variable in (\ref{gamma_xi}) is
\begin{eqnarray}
Y_V &:=& \frac 1{a(r_\epsilon) \sqrt{hl}} \sum_{(i,j)\in A_V \times B_V} \sum_k Z_{ij,k} \label{def-Y} \\
\mbox{{\rm with} } \;Z_{ij,k}&=& \log \cosh(u_k \cdot \frac{x_{ij,k}}{\epsilon \sigma_k}) -\frac{u_k^2}2 + \frac{u_k^4}4, \nonumber  
\end{eqnarray}
where $u_k = \theta_k^\star/(\ep \sigma_k)$ and $\{\theta_k^\star \}_k$ is the optimal sequence defined in Proposition \ref{sol-opt-contr}. Recall that $\theta_k^\star$ is null for $k>T$ and thus the sum over $k$ contains a finite number of non null elements. Moreover, due to (\ref{cond:a2}), recall that we  have \begin{eqnarray}
\sum_k u_k^4 = 2 a^2(r_\ep). \label{uk-opt}\end{eqnarray}

1. We shall prove that $\P_0(\Gamma^C) \to 0$. Let us write more conveniently
\begin{eqnarray*}
\Gamma^C &=& \bigcup_{\xi \in T_{M,N}(m,n)}   \bigcup_{\begin{array}{ll}\delta_1 m\leq h \leq m, \\
\delta_1 n\leq l \leq n \end{array}} 
\bigcup_{V \subset \xi} \{Y_V > T_{hl}\}\\
&=& \bigcup_{\begin{array}{ll}\delta_1 m\leq h \leq m, \\
\delta_1 n\leq l \leq n \end{array}} 
\bigcup_{V \in T_{M,N}(h,l)} \{Y_V > T_{hl}\}.
\end{eqnarray*}
Therefore, we have
\begin{eqnarray*}
\P_0(\Gamma^C)  &\leq & \sum_{\begin{array}{ll}\delta_1 m\leq h \leq m, \\
\delta_1 n\leq l \leq n \end{array}} \sum_{V \in T_{M,N}(h,l)} \P_0 (Y_V T_{hl}> T^2_{hl})\\
&\leq &\sum_{\begin{array}{ll}\delta_1 m\leq h \leq m, \\
\delta_1 n\leq l \leq n \end{array}} \sum_{V \in T_{M,N}(h,l)} e^{-T_{hl}^2} \E_0 (\exp(\sum_{(i,j)\in A_V \times B_V} 
\frac{T_{hl}}{a(r_\ep)\sqrt{hl}} \sum_k Z_{ij,k} ))\\
&\leq & \sum_{\begin{array}{ll}\delta_1 m\leq h \leq m, \\
\delta_1 n\leq l \leq n \end{array}} 
C_M^h C_N^l e^{-T_{hl}^2} \left(\E_0(\exp(\frac{T_{hl}}{a(r_\ep)\sqrt{hl}} \sum_k 
Z_{11,k}))\right)^{hl}.
\end{eqnarray*}

 Using Equation  (\ref{uk-opt}) and applying  Lemma~\ref{generating_cosh} for 
 $\lambda = T_{hl}/(a(r_\ep)\sqrt{hl})$ which is  $O(1)$,   one obtains 
% and by the definition of $T_{hl}$,   
\begin{eqnarray*}
\E_0(\exp(\frac{T_{hl}}{a(r_\ep)\sqrt{hl}} \sum_k Z_{11,k}))
%&=& \E_0(\frac{T_{hl}}{a(r_\ep)\sqrt{hl}} \sum_k  
%(\log \cosh(u_k \cdot \eta_{11,k} ) -\frac{u_k^2}2 + \frac{u_k^4}4))\\
%&=& \exp( \frac{T^2_{hl}}{4a^2(r_\ep)hl} \sum_k u_k^4 (1+o(1)))\\
&=& \exp( \frac{T^2_{hl}}{2hl} (1+O(a(r_\ep) \max_k w_k^\star)). 
\end{eqnarray*}
Recall that $T^2_{hl} = 2\log(C_M^h C_N^l)(1+A_\ep)+2\log(mn)$ where $A_\ep = o(1)$ 
%is such that $a(r_\ep) \max_k w_k^\star = o(A_\ep)$ 
by \nref{Aep}. Therefore 
\begin{eqnarray*}
\P_0(\Gamma^C)  &\leq & 
\sum_{\begin{array}{ll}\delta_1 m\leq h \leq m, \\
\delta_1 n\leq l \leq n \end{array}} 
C_M^h C_N^l e^{-T_{hl}^2} 
(\exp(\frac{T^2_{hl}}{2hl}(1+O(a(r_\epsilon )\max_k w_k^\star ) ) ))^{hl}\\
&= & \sum_{\begin{array}{ll}\delta_1 m\leq h \leq m, \\
\delta_1 n\leq l \leq n \end{array}} 
 \exp(- \frac{T_{hl}^2}2 +\frac{T_{hl}^2}2 O(a(r_\ep) \max_k w_k^\star)  + \log(C_M^h C_N^l))\\
&\leq & \sum_{\begin{array}{ll}\delta_1 m\leq h \leq m, \\
\delta_1 n\leq l \leq n \end{array}} 
\exp(-\log (C_M^h C_N^l )(A_\ep -a(r_\ep) \max_k w_k^\star )-\log(mn))
=o(1),
\end{eqnarray*}
for large enough $m,\, n,\, M$ and $N$, as we have both $A_\ep \cdot \log (C_M^h C_N^l ) \to \infty$ and $a(r_\ep) \max_k w_k^\star = o(A_{1,\ep}) = o(A_{\ep})$, by \nref{Aep}. 
\\

\noindent 2. We have
\begin{eqnarray*}
\E_0(\tilde L_{\overline \pi} (\overline X)) &=& \E_0(
\frac 1{C_M^m C_N^n} \sum_{\xi \in T_{M,N}(m,n)} \frac{d\P_\xi}{d\P_0}(\overline X) \1(\Gamma_\xi) )
= \P_\xi(\Gamma_\xi),
\end{eqnarray*}
which tends to 1 if and only if $\P_\xi(\Gamma_\xi^C) \to 0$. As we can write
$$
\P_\xi(\Gamma_\xi^C) \leq \sum_{\begin{array}{ll}\delta_1 m\leq h \leq m, \\
\delta_1 n\leq l \leq n \end{array}} \sum_{V \subset \xi} \P_V(Y_V >T_{hl}),
$$
where $\P_V$ is such that   
$$\frac{d\P_V}{d\P_0}( \overline X) = \prod_{(i,j) \in A_V \times B_V}  \prod_k \exp(-\frac{u_k^{2}}{2}) \cosh(u_k \frac{ x_{ij,k}}{\epsilon \sigma_k}) = \exp (Y_V a(r_\epsilon) \sqrt{hl} - lh   \frac{a^{2}(r_\epsilon)}{2}).$$ 
Then, applying  Lemma \ref{generating_cosh}, one obtains  for any positive $\lambda$ such that $(\lambda +1) = O(1)$, 
\begin{eqnarray}
\P_V(Y_V >T_{hl}) 
%& =& \E_0 ( \exp (Y_V a(r_\epsilon) \sqrt{hl} - lh   a^{2}(r_\epsilon) ) 1_{Y_V a(r_\epsilon) \sqrt{hl}>T_{hl} a(r_\epsilon) \sqrt{hl}}  ) \nonumber \\
& =& \P_V (Y_V a(r_\epsilon) \sqrt{hl}>T_{hl} a(r_\epsilon) \sqrt{hl}) \nonumber \\
& \leq & \E_V [ \exp( \lambda Y_V a(r_\epsilon) \sqrt{hl}) ]\exp  (-\lambda T_{hl} a(r_\epsilon) \sqrt{hl}) \nonumber \\
& =& \E_0 [ \exp( (\lambda +1)  Y_V a(r_\epsilon) \sqrt{hl}) ] \exp( -lh \frac{a^{2}(r_\epsilon)}{2}- \lambda T_{hl} a(r_\epsilon) \sqrt{hl}) \nonumber \\
& =& \exp ((\lambda +1)^{2}  lh \frac{a^{2}(r_\epsilon)}{2} (1 +o(1)) -lh \frac{a^{2}(r_\epsilon)}{2}- \lambda T_{hl} a(r_\epsilon) \sqrt{hl}) \label{measure-nu}
\end{eqnarray}
The minimum value for the right-hand side of (\ref{measure-nu}) is 
$$ 
\exp(-\frac{(T_{hl} - a(r_\epsilon) \sqrt{hl})^{2}}{2}(1 + o(1)))
$$
which is 
achieved for 
$\lambda =  \frac{T_{hl} }{a(r_\epsilon) \sqrt{hl}} -1$. Note that $\lambda$ satisfies $(\lambda +1)= O(1)$     
and that asymptotically $ \frac{T^{2}_{hl} }{a^{2}(r_\epsilon) hl} > \frac{2}{2 - \delta} > 1$.
%  which is required 
%to use Lemma \ref{generating_cosh}. 

%
%we have to study the distribution of $Y_V$ under $\xi$. 
%
%We see that $\E_\xi(Y_V) = a(r_\ep)\sqrt{hl}(1+o(1))$.
%We do the similar computations as for the previous question and get
%\begin{eqnarray*}
%\P_\xi(Y_V > T_{hl}) 
%&\leq & \P_\xi((Y_V -\E_\xi(Y_V))(T_{hl} - a(r_\ep)\sqrt{hl}) > (T_{hl}-a(r_\ep)\sqrt{hl})^2) \\
%&\leq & \exp( - (T_{hl} - a(r_\ep)\sqrt{hl} )^2) \E_\xi(e^{(T_{hl}-a(r_\ep)\sqrt{hl})(Y_V-\E_\xi(Y_V))})\\
%&\leq & \exp(- (T_{hl} - a(r_\ep)\sqrt{hl} )^2) \E^{hl}_{\boldsymbol{\theta}_{11}}(
%\exp(\frac{T_{hl}-a(r_\ep)\sqrt{hl}}{a(r_\ep)\sqrt{hl}} \sum_k (Z_{11,k}- \E(Z_{11,k})))).
%% \P_\xi(T_{hl}(Y_V - a(r_\ep)\sqrt{hl})>T^2_{hl}-T_{hl} )\\
%% &\leq & e^{-T^2_{hl}+T_{hl} a(r_\ep)\sqrt{hl}}\E_\xi( \exp(T_{hl}(T_1+T_2)))\\
%% &\leq & \exp(-T^2_{hl}+T_{hl} a(r_\ep)\sqrt{hl} + \frac{T^2_{hl}}2(1+o(1)) 
%% +2 T^2_{hl} a(r_\ep)(1+o(1)))\\
%% &\leq & \exp(-\frac{T^2_{hl}}2(1+o(1)).
%\end{eqnarray*}
%By Lemma~\ref{}, we have for $\lambda = (T_{hl}-a(r_\ep)\sqrt{hl})/a(r_\ep)\sqrt{hl} =O(1) $ that {\bf voila le resultat que j espere a rajouter dans le lemme 6.2 - sous quelles conditions?}
%$$
%\E_{\boldsymbol{\theta}_{11}}(
%\exp(\frac{T_{hl}-a(r_\ep)\sqrt{hl}}{a(r_\ep)\sqrt{hl}} \sum_k (Z_{11,k}- \E(Z_{11,k}))))
%\leq \exp(\frac{(T_{hl}-a(r_\ep)\sqrt{hl})^2}{2hl}(1+o(1))).
%$$ 
In conclusion,
\begin{eqnarray*}
\P_\xi(\Gamma_\xi^C) &\leq & 
\sum_{\begin{array}{ll}\delta_1 m\leq h \leq m, \\
\delta_1 n\leq l \leq n \end{array}}
C_m^h C_n^l \exp(-\frac 12 (T_{hl} - a(r_\ep)\sqrt{hl})^2(1+o(1)))\\
&\leq & \sum_{\begin{array}{ll}\delta_1 m\leq h \leq m, \\
\delta_1 n\leq l \leq n \end{array}}
C_m^h C_n^l \exp(-\frac {c(\delta)}2 T_{hl}^2(1+o(1))),
\end{eqnarray*}
for some $c(\delta)>0$ small with $\delta$.  Therefore, due to the asymptotic value of $T_{hl}$,  
% We use that $\log(C_m^h C_n^l) = O(m+n)$ for the values of $h$ and $l$ in the previous sum. 
 $\P_\xi(\Gamma_\xi^C) \leq  o(1)$.
\\

\noindent 3. We have, for $\xi_1= \1((i,j)\in A_1 \times B_1)$ and $\xi_2= \1((i,j) \in A_2 \times B_2)$,
\begin{eqnarray*}
\E_0(\tilde L_{\overline \pi}^2 (\overline X)) 
& = & \frac 1{(C_M^m C_N^n)^2} \sum_{\xi_1} \sum_{\xi_2}
g(h,l), \mbox{ where }\\
g(h,l)&=&\prod_{(i_1,j_1) \in A_1 \times B_1} \prod_{(i_2,j_2) \in A_2 \times B_2}
\prod_k \exp(-\frac{{\theta_k^{\star}}^2}{\epsilon^2 \sigma_k^2})\\
&& \cdot \E_0\left( \cosh (x_{i_1 j_1,k} \frac{\theta_k^{\star}}{\epsilon^2 \sigma_k^2}) 
\cosh (x_{i_2 j_2,k} \frac{\theta_k^{\star}}{\epsilon^2 \sigma_k^2}) \1(\Gamma_{\xi_1}\cap \Gamma_{\xi_2})\right)\\
&=& \prod_{(i,j) \in (A_1 \cap A_2) \times (B_1 \cap B_2) }
\prod_k \exp(-\frac{{\theta_k^{\star}}^2}{\epsilon^2 \sigma_k^2})
\E_0\left( \cosh^2 (x_{i j,k} \frac{\theta_k^{\star}}{\epsilon^2 \sigma_k^2}) 
\1(\Gamma_{\xi_1}\cap \Gamma_{\xi_2})\right)
\end{eqnarray*}
and the function $g$ depends on the sets $A_{1},A_{2}$ and $B_{1},B_{2}$ only through the
number $h$ of common rows of $A_1$ and $A_2$ and the number $l$ of common columns of $B_1$ and $B_2$.
After some combinatorics we can write
\begin{eqnarray*}
\E_0(\tilde L_{\overline \pi}^2 (\overline X)) 
& = & \E(g(H,L)),
\end{eqnarray*}
where $H$ and $L$ are independent random variables having hypergeometric distribution $\mathcal{HG}(M,m,m)$ and $\mathcal{HG}(N,n,n)$, respectively. Let us see that, for any $0 \leq h \leq m$ and $0\leq l \leq n$,
\begin{eqnarray*}
\log (g(h,l)) &\leq & \sum_{(i,j) \in (A_1 \cap A_2) \times (B_1 \cap B_2) }
\log \left(\prod_k \exp(-\frac{{\theta_k^\star}^2}{\epsilon^2 \sigma_k^2})
\E_0\left( \cosh^2 (x_{i j,k} \frac{\theta_k^\star}{\epsilon^2 \sigma_k^2}) \right)
\right)\\
& =& h l \log \left(\prod_k \exp(-\frac{{\theta_k^{\star}}^2}{\epsilon^2 \sigma_k^2})
\frac{1}{2} \left( \exp(  \frac{2{\theta_k^{\star}}^2}{\epsilon^2 \sigma_k^2}) +1 \right)
\right) \\
& =& h l \log \left(\prod_k \cosh ( \frac{{\theta_k^{\star}}^2}{\epsilon^2 \sigma_k^2})
\right)
= : hl \cdot D.
\end{eqnarray*}
Therefore, $\E(g(H,L)) \leq \E(e^{HL \cdot D})$ for $D$ which has the following asymptotic equivalent
\begin{eqnarray}
D & =& \sum_k\log \left( 1 +2\sinh^{2} ( \frac{{\theta_k^{\star}}^2}{2\epsilon^2 \sigma_k^2})
\right) 
 = \sum_k\log \left( 1 +2 (\frac{{\theta_k^{\star}}^2}{2\epsilon^2 \sigma_k^2})^{2}(1+ o(1)) 
\right)  \nonumber \\
& =& \sum_k (\frac{{\theta_k^{\star}}^4}{2\epsilon^4 \sigma_k^4})(1+ o(1))
 = \frac{V_\ep}{\ep^{2}} \sum_k w_k^{\star} \frac{(\theta_k^{\star})^2}{\epsilon^2 \sigma_k^2}(1+ o(1)) \nonumber\\
&=&  \frac{V_\ep^{2}}{\ep^{4}}(1+ o(1))= a^{2}(r_\ep) (1+ o(1)),  \label{D-eval}
%\\
%&=& - \sum_k \frac{{\theta_k^{\star}}^2}{\epsilon^2 \sigma_k^2}
%+ \sum_k \log \left( E_0\left( \cosh^2 (x_{i j,k} \frac{\theta_k^*}{\epsilon^2 \sigma_k^2}) \right) \right)\\
%&=& -a(r_\epsilon) + \sum_k \log \left( \frac 12 +\frac 12 \cosh (\frac{2 \theta_k^2}{\epsilon^2 \sigma_k^2})\right).
%\end{eqnarray*}
%As we have $\theta_k^2/(\epsilon^2 \sigma_k^2) = w_k a(r_\epsilon)$ is an $o(1)$, we use asymptotic equivalents to get
%\begin{eqnarray*}
%A &=& -a(r_\epsilon) + \sum_k \log \left( 1+\frac 12 \frac {(4w_k a(r_\epsilon))^2}2(1+o(1)) \right)\\
%&=& -a(r_\epsilon) + \sum_k 4 w_k^2 a^2(r_\epsilon) (1+o(1))\\
%&=& -a(r_\epsilon) + 2 a^2(r_\epsilon) (1+o(1)).
\end{eqnarray}
which holds under Assumption (\ref{Cond-a-m-n}). %  $\displaystyle{\max_k w_k^{\star}a(r_\epsilon)} =o(1)$.
Indeed, this Assumption implies that
$$
\displaystyle{\max_k w_k^{\star}a(r_\epsilon)} \leq \frac{v\sigma_T^2}{2\ep^2} \leq C \ep^{-2} r_\ep^{2+(2s+1)/\tau}
\leq C \ep^{-2} \tilde{r_\ep}^{2+(2s+1)/\tau} =o(1).
$$

We shall split $\E(g(H,L))$ into the sum $I_1+I_2$, where
\begin{eqnarray*}
I_1 &=& \E(g(H,L)\cdot \1(HD <1)),\\
I_2 &=& \E(g(H,L) \cdot \1(HD \geq 1)).
\end{eqnarray*} 
For $I_1$, we use the stochastic ordering in Lemma~\ref{lem:stoch-order} and replace the hypergeometric distributions of $H$ and $L$ with binomial distributions of $\tilde H \sim Bin(m,\tilde p)$ and $\tilde L \sim Bin(n,\tilde q)$, respectively. We have
\begin{eqnarray*}
I_1 &\leq & \E(\E(e^{\tilde L \tilde H D} \1(\tilde H D < 1) \vert \tilde H)\\
&\leq & \E((1+\tilde q (e^{\tilde H D} -1))^n \1(\tilde H D < 1))\\
&\leq & \E(\exp(C n\tilde q \tilde H D )) = (1 + \tilde p (e^{C n\tilde q D} -1))^m,
\end{eqnarray*}
for some constant $C>0$. By assumption (\ref{lb2}),  (\ref{aasymp}) and (\ref{D-eval}), $Dn \asymp \log(p^{-1})$, which implies that $D n \tilde q \asymp (\tilde q \log(p^{-1})) $ and this is an $o(1)$ by assumption (\ref{lb1}). 
So, by assumption (\ref{lobo1}), $I_1 \leq \exp(C mn\tilde p \tilde q D) = 1+o(1)$. 

The rest of the Section is devoted to the proof that $I_2 = o(1)$. We shall further split the expected value into the sum of $I_{21}+I_{22}$, where
\begin{eqnarray*}
I_{21} &=& \E(g(H,L) \cdot \1(HD \geq 1) \cdot \1(L< n \delta_1)),\\
I_{22} &=& \E(g(H,L) \cdot \1(HD \geq 1) \cdot \1(L \geq n \delta_1)), 
\end{eqnarray*}
for some fixed $\delta_1>0$, small enough such that $Dn \delta_1 \leq \log(p^{-1})/2$ and that $Dm\delta_1 \leq \log(q^{-1})/2$.
On the one hand
\begin{eqnarray*}
I_{21}&\leq & \sum_{D^{-1} <h\leq m,\, 0\leq l<n \delta_1} e^{hlD} \P_{\mathcal{HG}(M,m,m)}(H = h) \P_{\mathcal{HG}(N,n,n)}(L = l)\\
&\leq & \sum_{D^{-1} <h\leq m,\, 0\leq l<n \delta_1} e^{h(lD - \log(p^{-1})(1+o(1)))} ,
\end{eqnarray*}
as $\P_{\mathcal{HG}(N,n,n)}(L = l) \leq 1$ and by using Lemma 
5.3 in \cite{BI} for $\log (\P_{\mathcal{HG}(M,m,m)}(H = h)) \leq h \log(p) (1+o(1))$. Now, under the constraints in the sum, $lD \leq Dn\delta_1 \leq (1/2 +o(1))\log(p^{-1})$. This implies that
$h(lD - \log(p^{-1})(1+o(1))) \leq -h \log(p^{-1}) (1/2+o(1)) \leq -D^{-1} \log(p^{-1}) (1/2+o(1)) \asymp -n$. Therefore,
$$
I_{21}\leq mne^{-Bsn}, \mbox{{\rm  for some fixed} } B>0,
$$
and this is an $o(1)$.

We can also split $I_{22}$ into the sum of $I_{221}+I_{222}$, where 
\begin{eqnarray*}
I_{221} &=& \E(g(H,L) \cdot \1(HD \geq 1, H <m\delta_1) \cdot \1(L \geq n \delta_1)),\\
I_{222} &=& \E(g(H,L) \cdot \1(HD \geq 1, H \geq m\delta_1) \cdot \1(L \geq n \delta_1)). 
\end{eqnarray*}
It is easy to check that $I_{221} = o(1)$ as we previously did for $I_{21}$.

On the other hand, we can write
$$
I_{222} = \E(e^{HL D} \1(\mathcal{H})), \mbox{ where } \mathcal{H} = \{(h,l),\,  m \delta_1 \leq h \leq m, n \delta_1 \leq l \leq n\}.
$$
Note that under the event $\mathcal{H}$ we have $T^2_{hl} = (2+\delta) (h/m \cdot m \log(p^{-1})+ l/n \cdot n \log(q^{-1})) \geq \delta_1 T^2_{mn}$.

We divide again the set $\mathcal{H}$ in disjoint sets 
$$
\mathcal{H}_1=\{(h,l)\in \mathcal{H}: T^2_{hl} > 2 T^2_{mn} \frac{hl}{mn}\}
\mbox{ and } 
\mathcal{H}_2=\{(h,l)\in \mathcal{H}: T^2_{hl} \leq 2 T^2_{mn} \frac{hl}{mn}\}.
$$

Let us go back to $L_{\overline \pi}(\overline X)$ and write it as
\begin{eqnarray*}
L_{\overline \pi}(\overline X) &=& \frac 1{C_M^m C_N^n} \sum_{\xi \in T_{M,N}(m,n)} 
\exp\left(  \sum_{(i,j) \in A_\xi \times B_\xi } \sum_k (-\frac{{\theta_k^{\star}}^2}{2\epsilon^2 \sigma_k^2} 
+\log \cosh (x_{i j,k} \frac{\theta_k^{\star}}{\epsilon^2 \sigma_k^2})) \right)\\
&=&\frac 1{C_M^m C_N^n} \sum_{\xi \in T_{M,N}(m,n)} 
\exp\left(  \sum_{(i,j) \in A_\xi \times B_\xi} \sum_k 
(\log \cosh (\frac{x_{i j,k}}{\epsilon \sigma_k} u_k) -\frac{u_k^2}{2} )
 \right)\\
&=&\frac {\exp(-a^2(r_\epsilon ) mn/2 )}{C_M^m C_N^n} \sum_{\xi \in T_{M,N}(m,n)} 
\exp\left( \sum_{(i,j) \in A_\xi \times B_\xi} \sum_k 
(\log \cosh (\frac{x_{i j,k}}{\epsilon \sigma_k} u_k) -\frac{u_k^2}{2} +
\frac{u_k^4}4 )\right),
\end{eqnarray*}
where $u_k = \theta_k^{\star}/(\epsilon \sigma_k )$ is such that $\sum_k u_k^4 = 2 a^2(r_\epsilon)$ (see (\ref{cond:a2})).

Now, we give a tighter upper bound for $g(h,l)$ than the one used for $I_1$. 
Using the same notation as to define  $Y_V$ in (\ref{def-Y}) and for any matrix $\xi$, we  define the random variable $Y_\xi= \frac 1{a(r_\epsilon) \sqrt{hl}} \sum_{(i,j)\in A_\xi \times B_\xi} \sum_k Z_{ij,k}$.  
  Then,
we write
\begin{eqnarray}
%\E_0 (\tilde L_{\overline \pi}(\overline X)^2 I(\Gamma))
g(h,l)
&=& e^{-a^2(r_\ep) mn }
\E_0(e^ {a(r_\ep) \sqrt{mn} (Y_{\xi_1}+Y_{\xi_2})} \1(\Gamma_{\xi_1} \cap \Gamma_{\xi_2} ) )\nonumber \\
&\leq & e^{-a^2(r_\ep) mn }\E_0(e^ {a(r_\ep) \sqrt{mn} (Y_{\xi_1}+Y_{\xi_2})}
\1(Y_{\xi_1} \leq T_{mn}, Y_{\xi_2} \leq T_{mn}) )\nonumber \\
&\leq & e^{-a^2(r_\ep) mn +2 T_{mn} J} \E_0(e^ {(a(r_\ep) \sqrt{mn}- J) (Y_{\xi_1}+Y_{\xi_2})}),\label{second_bound}
\end{eqnarray}
for some $J>0$ that we will choose later on. In order to deal with $I_{222}$, we keep in mind that we consider only submatrices $\xi_1$ and $\xi_2 $ having $h$ common rows and $l$ common columns, such that $(h,l) \in \mathcal{H}$.
Denote by $V$ the submatrix of common rows and columns for $\xi_1$ and $\xi_2$,
that is
$$
V=\I((i,j) \in (A_{\xi_1}\times B_{\xi_1}) \cap (A_{\xi_2}\times B_{\xi_2})),
$$ 
and by $V_1 = \xi_1-V$ (respectively $V_2=\xi_2-V$). Therefore, 
$$
\sqrt{mn} (Y_{\xi_1}+Y_{\xi_2}) = \sqrt{mn-hl}(Y_{V_1}+Y_{V_2})+2\sqrt{hl} Y_V.
$$ Replace this into the equation (\ref{second_bound}) and get by Lemma~\ref{generating_cosh}
\begin{eqnarray*} 
&& I_{222} \\
&\leq & \sum_{(h,l)\in \mathcal{H}}\exp(-a^2(r_\ep) mn +2 T_{mn} J + (a(r_\ep)\sqrt{mn}-J)^2 (1+\frac{hl}{mn}))\P_{\mathcal{HG} (M,m,m)}(h)\P_{\mathcal{HG}(N,n,n)}(l)\\
&\leq & \sum_{(h,l)\in \mathcal{H}}\exp(-a^2(r_\ep) mn +2 T_{mn} a(r_\ep ) \sqrt{mn} - \frac{T^2_{mn}}{1+hl/(mn)} - (h \log(p^{-1}) + l \log(q^{-1}))(1+o(1)))\\
&\leq & \sum_{(h,l)\in \mathcal{H}}\exp(- (a(r_\ep) \sqrt{mn} -T_{mn})^2 + \frac{T^2_{mn}hl}{mn + hl} - \frac{T^2_{hl}}{2}(1+o(1))),
\end{eqnarray*}
for $a(r_\ep)\sqrt{mn}-J = T_{mn}/(1+ hl/(mn))$.
Note that, for $\delta>0$ small enough, there exists $\delta_2>0$ such that 
$(a(r_\ep) \sqrt{mn} -T_{mn})^2 \geq \delta_2 T^2_{mn}$. Moreover, for $(h,l) \in \mathcal{H}_1$, 
\begin{eqnarray*}
\frac{T^2_{mn}hl}{mn + hl} - (h \log(p^{-1})+ l \log(q^{-1}))(1+o(1))
&\leq & \frac{T^2_{mn} hl}{mn+hl} -\frac{T^2_{hl}(1+o(1))}{2} \\ 
&\leq & T^2_{mn} \frac {hl}{mn}(\frac 1{1+ hl/(mn)} - 1)+o(T^2_{mn}),
\end{eqnarray*}
which is asymptotically less than $o(T^2_{mn})$. This implies that $I_{222}=o(1)$ over the set $\mathcal{H}_1$.

Finally, we give a yet slightly different upper bound for $g(h,l)$ in order to deal with $I_{222}$ when $(h,l)$ belongs to $\mathcal{H}_2$.
\begin{eqnarray*}
g(h,l)
&=& e^{-a^2(r_\ep) mn }
\E_0(e^ {a(r_\ep) \sqrt{mn} (Y_{\xi_1}+Y_{\xi_2})} \1(Y_{V} \leq T_{hl}) )\\
&\leq & e^{-a^2(r_\ep) hl }\E_0(e^ {2 a(r_\ep) \sqrt{hl} Y_{V}}
\1(Y_{V} \leq T_{hl}) )\\
&\leq & e^{-a^2(r_\ep) hl +T_{hl} J}\E_0(e^ {(2 a(r_\ep) \sqrt{hl}-J) Y_{V}
+J(Y_V - T_{hl})}
\1(Y_{V} \leq T_{hl}) )\\
&\leq &e^{-a^2(r_\ep) hl +T_{hl} J + (2 a(r_\ep) \sqrt{hl}-J)^2/2}.
\end{eqnarray*}
Take $J = 2 a(r_\ep) \sqrt{hl} - T_{hl}$ which is indeed positive for $(h,l)$ in $\mathcal{H}_2$ and obtain
$$
g(h,l) \leq \exp(- a^2(r_\ep) hl +2 a(r_\ep) \sqrt{hl} T_{hl}  - \frac{T^2_{hl}}{2} ).
$$
Moreover, denote $D^2_{hl} = h \log(p^{-1}) + l \log(q^{-1})$ and see that $D^2_{mn} hl/(mn) \leq D^2_{hl}$. 
%Therefore, $a^2(r_\ep )hl = (2-\delta) D^2_{mn} hl/(mn)$, $T^2_{hl} = (2+\delta) D^2_{hl}$ and $D^2_{hl} \leq (2+\delta) D^2_{mn} hl/(mn)$ under $\mathcal{H}_2$. 
We get
\begin{eqnarray*}
&& I_{222} \\
&\leq & \sum_{(h,l)\in \mathcal{H}}\exp( - a^2(r_\ep) hl +2 a(r_\ep) \sqrt{hl} T_{hl}  - \frac{T^2_{hl}}{2}
 )\P_{\mathcal{HG}(M,m,m)}(h)\P_{\mathcal{HG}(N,n,n)}(l)\\
&\leq & \sum_{(h,l)\in \mathcal{H}}\exp(- (a(r_\ep)\sqrt hl - T_{hl})^2 +\frac{T^2_{hl}}{2}- D^2_{hl} (1+o(1)))\\
&\leq & \sum_{(h,l)\in \mathcal{H}}\exp(- \frac{\delta^2}8 D^2_{hl} + o(1)D^2_{hl}
)=o(1).
\end{eqnarray*}

\hfill $\Box$

\subsection{Proof of Lemma \ref{generatrice-function}} 
For the sake of simplicity, we omit in this part the indices $i$ and $j$ so that $t_{ij,w^{\star}}$
%,$\boldsymbol{\theta}_{ij}$ 
and $\eta_{ij,k}$ are denoted by $t_{w^{\star}}$
%, $\boldsymbol{\theta}$ 
and $\eta_k$, respectively. 

%Set $u_k= \displaystyle{\frac{\theta_k}{\epsilon \sigma_k}}$. Note also that
%for $(i,j)$ such that $\xi_{ij}=0$, one may consider that $\boldsymbol{\theta}_{ij} \equiv 0$ which implies that $u_k =0, \forall k$. 
Under $H_0$,  observe that 
$t_{w^{\star}}= \sum_k w_k^{\star} (  \eta_{k}^{2} -1)$, with $\eta_k \stackrel{iid}{\sim} {\cal N} (0,1)$.  
Using the fact that  $\displaystyle{\mathbb{E}(e^{t \eta_k^2}) = 
 \frac{1}{(1-2t)^{1/2}}, 
\; \mbox{ for } t \leq \frac 1 2},
$
% 
%%Let us start with the distribution $\mathbb{P}_{\boldsymbol{\theta}}$. 
%%Note that $(u_k+\eta_k)^2$ has a noncentral $\chi^{2}(1)$ distribution with noncentrality parameter $u_{k}^{2}$.
% Therefore, it has the following moment generating function
%$$
%\mathbb{E}(e^{t(u_k+\eta_k)^2}) = 
%\exp \left( \frac{t u_k^{2}}{1-2t} \right) \frac{1}{(1-2t)^{1/2}}, 
%\quad \mbox{ for } t \leq \frac 1 2.
%$$
we obtain for $\lambda$ such that 
$\displaystyle{\lambda \max_k w_k^{\star}} =o(1)$,  
\begin{eqnarray}
\mathbb{E}_{0} (\exp (\lambda t_{w^{\star}})) 
%&= &\mathbb{E}_{\boldsymbol{\theta}} (\prod_{k \in \mathbb{Z}} \exp(-\lambda w_k^\star +\lambda w_k^\star \eta_k^2  )) \nonumber \\
&=&\prod_{k \in \mathbb{Z}} \exp(-\lambda  w_k^{\star} -\frac 12 \log(1-2\lambda w_k^\star)) \nonumber \\
%& =& \prod_{k \in \mathbb{Z}} \exp \left(  \lambda^{2}(w_k^{\star})^{2} +o(\lambda^{2}(w_k^{\star})^{2}) \right)  \nonumber \\
& =& \exp \left( 
  \lambda^{2} \sum_k (w_k^{\star})^{2} (1 +o(1)) \right) \nonumber \\ 
  & =&  \exp \left( 
  \frac{\lambda^{2}}{2} (1 +o(1)) \right), \label{terme-1} 
\end{eqnarray}
where the last equality holds since  $\displaystyle{\sum_k (w_k^{\star})^{2} =\frac{1}{2}}$.  This ends the proof. \endproof

\subsection{Proof of Proposition~\ref{sol-opt-contr}}\label{sec:Lemma2.1}
These computations can be found in Ingster and Suslina~\cite{IS.02a}, but we give the sketch of proof for the convenience of the reader.

Let us change variables in problem (\ref{max-min_pb}), by defining  $v_k = \displaystyle{\frac{\theta_k^2}{\sigma_k^2\sqrt{2}}}$,
% $u_k = w_k\sqrt{2}$, 
for all $k \in \mathbb{Z}$. We have $\{\theta_k\}_k$ belongs to $\Sigma(\tau,r_\epsilon)$ if and only if $\{v_k\}_k$ belongs to $\tilde \Sigma (\tau , r_\epsilon)$, where
\begin{equation*}\label{V+}
{\tilde \Sigma }(\tau,r_\ep) %=\{v: q\in \Theta(\tau, r_\e)\}=
=\left\{\{v_k\}_k \in l_1(\mathbb{Z}): v_k\geq 0;\; (2 \pi)^{2 \tau} \sum_{k\in
\mathbb{Z}}\vert k \vert^{2 \tau} \sigma_k^2 v_k \leq \frac 1{\sqrt{2}}; \;
  \sum_{k\in \mathbb{Z}} v_k \sigma_k^2 \geq \frac{r_\epsilon^2}{\sqrt{2}} \right\}.
\end{equation*}
The problem (\ref{max-min_pb}) is equivalent to
\begin{eqnarray*}
&&\frac{\sqrt{2}}{\ep^{2}}\sup_{\{w_k\}_k : \sum_k w_k^2=1/2, w_k\ge 0}\quad \inf_{\{v_k\}_k \in {\tilde \Sigma }(\tau,r_\ep) }\sum_k w_k v_k \\
&=&\frac{\sqrt{2}}{\ep^{2}} \sup_{\{w_k\}_k : \sum_k w_k^2\le
1/2, w_k\ge 0}\quad   \inf_{\{v_k\}_k \in {\tilde \Sigma }(\tau,r_\ep) }\sum_k w_k v_k \\
& =& \frac{\sqrt{2}}{\ep^{2}} \inf_{\{v_k\}_k\in {\tilde \Sigma }(\tau,r_\ep) } \quad \sup_{\{w_k\}_k : \sum_k w_k^2\le
1/2, w_k\ge 0} \sum_k w_k v_k ,
%& =& \frac{1}{\ep^{2}}\inf_{\{v_k \}_k \in {\tilde \Sigma }(\tau,r_\ep) }(\stackrel{\Delta}{=}\frac{1}{\ep^{2}} V_\ep.
%&=&\frac{1}{\sqrt{2}\;\ep^2}\inf_{\theta\in \Sigma (\tau,r_\ep)}(\sum_k \theta_k^4)^{1/2}=a(r_\e).
\end{eqnarray*}
by the minimax theorem on convex sets. Now, use the Cauchy-Schwarz inequality to see that 
$$
 \sqrt{2}\sup_{\{ w_k \}_k :\sum_k w_k^2\leq 1/2, w_k\geq 0}\sum_k w_k v_k = (\sum_k
v_k^2)^{1/2}
$$
and the equality holds for $w_k = v_k (2 \sum_k
v_k^2)^{-1/2}$. As we denoted by $V_\epsilon^2 = \sum_k (v_k^\star)^2$
we get $ w_k^\star = v_k^\star/(\sqrt{2} V_\epsilon) $, which is equivalent to
$$ 
w_k^\star = \frac{(\theta^\star_k)^2}{2 \sigma_k^2 \; V_\ep} , \mbox{ for all } k \in \mathbb{Z}.
$$
It follows that solving the 
%The last step is to compute $v_k^\star$ and deduce $\theta_k^\star$, $w_k^\star$ and an asymptotic equivalent of $V_\epsilon$.
%  Observe that  
%  \begin{eqnarray}
% \inf_{ \{\theta_k\}_{k \in \mathbb{Z}} \in
%   \Sigma (\tau,r_{\epsilon})} \sum_{k \in \mathbb{Z}} w_k^\star \; \;\frac{\theta^2_{k}}{\ep^{2} \sigma_k^2}& \geq & V_\epsilon\frac{1}{\ep^{2}} . \label{PB_minim} 
% \end{eqnarray}
% Set $a(r_\ep)=\displaystyle{  V_\epsilon\frac{1}{\ep^{2}}}$. 
 problem (\ref{max-min_pb}) reduces to solve the optimization program
 \begin{equation*}
 \inf_{\{v_k\}_k, v_k \geq 0} \sum_k v_k^2 +\lambda_1 \left( 
\sum_k (2\pi |k|)^{2\tau} \sigma_k^2 v_k -\frac 1{\sqrt{2}}\right)
-\lambda_2 \left( \sum_k \sigma_k^2 v_k - \frac{r_\epsilon^{2}}{\sqrt{2}}\right)
 \end{equation*}
By the  Lagrangian multipliers rules, one gets for $\lambda_1 \in \mathbb{R}$ and $\lambda_2 \in \mathbb{R}$ the following system of equations
$$
\left\{ \begin{array}{l} 
2 \sum_k v_k + \lambda_1\sqrt{2} (2 \pi)^{2 \tau}(\sum_k \vert k\vert^{2 \tau} \sigma_k^{2} ) - \lambda_2 \sqrt{2} \sum_k \sigma_k^{2} =0 \\
\sqrt{2}(2 \pi)^{2 \tau}  \sum_k v_k \vert k\vert^{2 \tau} \sigma_k^{2}  = 1 \\
\sqrt{2} \sum_k v_k \sigma_k^{2} =r_\ep^2
\end{array} \right. 
$$  

Put, for all $k \in \mathbb{Z}$,  $\displaystyle{v_k = v \sigma_k^{2} \left(1 - (\frac{\vert k \vert}{T})^{2 \tau}\right)_+}$, 
where $v = \frac{\lambda_2}{\sqrt{2}}$, 
$T =\frac 1{2 \pi} \left(\frac{\lambda_2}{\lambda_1}\right)^{1/(2\tau)}$ and $(x)_+=\max(0,x)$. 

We evaluate the solution of the previous system as $T$ goes to infinity. Using $\sigma_k \sim \vert k\vert^s$ for $\vert k\vert$ large enough and some $s>0$,  the last two equations in the previous system become
$$
\left\{ \begin{array}{l} 
\kappa_2 v T^{2\tau + 4 s +1} \sim 1\\
\kappa_1 v T^{4s+1} \sim r_\epsilon^2,
\end{array} \right. 
$$
that  gives 
$$
T \sim \left(\frac{\kappa_1}{\kappa_2} \right)^{\frac 1{2\tau }} r_\ep^{- \frac 1{\tau}} \mbox{ and }
v \sim \frac{1}{\kappa_1}  \left(\frac{\kappa_2}{\kappa_1}\right)^{\frac{4s+1}{2\tau }} r_\ep ^{2 + \frac{4s+1}\tau}.
$$ 
Note that $T\to \infty$ provided that $r_\ep \to 0$. 
It further gives
\begin{eqnarray*}
V_\ep^2 &=& \sum_{|k|\leq T} (v_k^\star)^2 \sim 2 v^2 T^{4s+1} \kappa_3
\sim  c(\tau,s)^2 r_\ep ^{4+\frac{4s+1}{ \tau}}.
\end{eqnarray*}
Finally, it is straightforward that 
 \begin{eqnarray*}
 \max_k w_k^{\star} & \leq & \frac{v}{V_\epsilon}  \max_{0\leq \vert k \vert \leq T}  \sigma_k^{2} \left(1 - (\frac{\vert k \vert}{T})^{2 \tau}\right)\nonumber\\ 
 %& < & \frac{v}{V_\epsilon}    \max_{0\leq \vert k \vert \leq T}  \sigma_k^{2}  \nonumber \\
 & \leq & \frac{v}{V_\epsilon}   \sigma_T^{2} 
 \asymp r^{\frac{4s+2 \tau +1}{\tau} - \frac{4s+4 \tau +1}{2\tau}-\frac{2s}{\tau}}_\ep 
 = r_\ep^{\frac 1 {2 \tau }} \stackrel{r_\ep \rightarrow 0}{\longrightarrow} 0. \end{eqnarray*}
\hfill $\Box$

%
%\bibliographystyle{plain}
%
%
%\bibliography{BibGhis}

\begin{thebibliography}{99}

\bibitem{BI} \textsc{Butucea, C.} and \textsc{Ingster, Yu.I.} (2012).
Detection of a sparse submatrix of a high-dimensional noisy matrix. 
 \textit{Bernoulli}, to appear.

\bibitem{CGPT} \textsc{Cavalier, L., Golubev, G.K., Picard, D. and Tsybakov, A.B.} (2002).
{Oracle inequalities for inverse problems. Dedicated to the memory of Lucien Le Cam.
\textit Ann. Statist.},  {\bf 30}, 843--874.

\bibitem{GI} \textsc{Gayraud, G.} and \textsc{Ingster, Yu.I.}  (2012).
Detection of sparse additive variable functions. 
\textit{Electronic Journal of Statistics}, {\bf 6}, 1409--1448.


\bibitem{IS.02a} \textsc{Ingster, Yu.I.} and \textsc{Suslina, I.A.} (2003).
\textit{Nonparametric goodness-of-fit testing under gaussian
models}. Lectures Notes in Statistics.  vol. 169., Springer-Verlag, New York.

\bibitem{KM.10}  \textsc{Klenke, A.}  and \textsc{Mattner, L.} (2010). Stochastic ordering of classical discrete distributions. {\it Adv. in Appl. Probab.},  {\bf 42}, 392--410. 

\bibitem{SN}
\textsc{Sun, X.} and \textsc{Nobel, A.B.} (2010). On the maximal size of
Large-Average and ANOVA-fit Submatrices in a Gaussian Random Matrix.
\textit{ArXiv: 1009.0562v1}

\bibitem{Shabalinetal}
\textsc{Shabalin, A.A.}, \textsc{Weigman, V.J.}, \textsc{Perou, C.M.} and \textsc{Nobel, A.B.} (2009).
Finding Large Average Submatrices in High Dimensional Data.
\textit{Annals of Applied Statistics}, {\bf 3}, 985--1012.


\bibitem{St.85}
\textsc{Stone, Ch.} (1985). Additive regression and other nonparametric
models. {\it Ann. Statist.},  {\bf 13},  689--705.

\end{thebibliography}
%
%
%%
%%\begin{thebibliography}{9}
%%\bibitem{IPT} \textsc{Ingster, Yu.I., Suslina, I.A.}   (2008), \textit{Nonparametric hypothesis testing for small type I errors.}
%%
%%\end{thebibliography}

\end{document}